\documentclass[a4paper,reqno]{amsart}
\newcommand{\ie}{\textit{i.e.} }
\newcommand{\cf}{\textit{cf.} }
\newcommand{\xdif}[1]{\mathrm{d}#1}
\newcommand{\Rref}[1]{\textup{\ref{#1}}}

\usepackage{amsmath}
\usepackage{amssymb}
\usepackage{amsfonts}

\usepackage[T1]{fontenc}
\usepackage[latin1]{inputenc}

\DeclareMathOperator\rank{Rank}
\def\XX{{X}}
\def\FF{{Y}}
\def\EE{{E}}

\newcommand{\ordreplat}{j}
\newcommand{\NN}{\mathbb{N}}
\newcommand{\RR}{\mathbb{R}}
\newcommand{\partialx}[1]{\frac{\partial}{\partial #1}}
\newcommand{\partialxy}[2]{\frac{\partial #1}{\partial #2}}

\newcommand{\sysedp}{\mathcal{E}^{\gamma,\delta}_{k,\ell}}
\newcommand{\sysedpl}[4]{\mathcal{E}^{#1,#2}_{#3,#4}}

\newcommand{\ouvedp}{O}

\newtheorem{thrm}{Theorem}[section]
\newtheorem{lmm}[thrm]{Lemma}
\newtheorem{prpstn}[thrm]{Proposition}

\theoremstyle{definition}
\newtheorem{dfntn}[thrm]{Definition}
\newtheorem{rmrk}[thrm]{Remark}
\newtheorem{xmpl}[thrm]{Example}
\newtheorem{cnjctr}[thrm]{Conjecture}

\begin{document}
{\small

\noindent\textsl{Manuscript submitted for publication} in :\\
\textbf{ESAIM: Control, Optimisation and Calculus of Variations}\\
URL: \texttt{http://www.edpsciences.org/cocv/}}

\bigskip\bigskip

\title[Monge parameterizations with 3 state and 2 controls]{
Flatness and Monge parameterization\\of two-input systems,
\\ control-affine with 4 states \\or general with 3 states}

\author{David Avanessoff} 
\author{Jean-Baptiste Pomet}

\dedicatory{INRIA, B.P. 93, 06902 Sophia Antipolis cedex, France\\
\textup{\texttt{David.Avanessoff@sophia.inria.fr}, 
\texttt{Jean-Baptiste.Pomet@sophia.inria.fr}}}

\date{April 21, 2005, revised November 30, 2005.} 

\begin{abstract}
  This paper studies Monge parameterization, or differential flatness,
  of  control-affine systems with four states and two
  controls. Some of them are known to be flat, and this implies admitting a
  Monge parameterization. 
Focusing on systems outside this class, we describe the only possible
structure of such a parameterization for these systems, and give a lower bound
on the order of this parameterization, if it exists. 
This lower-bound is good enough to recover the known results about
``$(x,u)$-flatness'' of these systems, with much more elementary techniques.
on the order of this parameterization, if it exists.

\medskip

\paragraph{\textsc{R{\'e}sum{\'e}}} On s'int{\'e}resse aux param{\'e}trisations de Monge, ou {\`a}
  la platitude, des syst{\`e}mes affines {\`a} quatre {\'e}tats et deux
  entr{\'e}es. Des travaux antérieurs caractérisent ceux de ces systèmes qui
  sont ``$(x,u)$-plats'', mais on ne sait pas si certains des sytèmes restants
  sont plats, ou non. La conjecture est qu'aucun n'est plat, ni
  Monge-paramétrable. Pour ces sytèmes, on montre que toute
  param{\'e}trisation est d'un type particulier, et on donne une borne
  inf{\'e}rieure sur l'ordre de cette param{\'e}trisation, suffisante pour
  retrouver, de mani{\`e}re beaucoup plus {\'e}l{\'e}mentaire, le résultat
  connu sur la ``$(x,u)$-platitude''.
\end{abstract}

\maketitle

\section{Introduction}

In control theory, after a line of research on exact linearization by dynamic
state feedback \cite{Isid-Moo-deL86,Char-Lev-Mar89,Char-Lev-Mar91}, the concept of
differential flatness was introduced in 1992 in\cite{Flie-Lev-Mar-R92cras} (see also
\cite{Flie-Lev-Mar-R95ijc,Flie-Lev-Mar-R99geo}).
Flatness is equivalent to exact linearization by dynamic state feedback of a
special type, called ``endogenous''~\cite{Flie-Lev-Mar-R92cras}, but, as
pointed out in that reference, it has its own interest, maybe more important than linearity.
An interpretation and framework for that notion is also proposed in
\cite{Aran-Moo-Pom95vars,Pome95vars,VanN-Rat-Mur98}; see
\cite{Mart-Rou-Mur02tri} for a recent review.

The \emph{Monge problem} (see the the survey article \cite{Zerv32}, published in 1932, that mentions
the prominent contributions \cite{Hilb12} and \cite{Cart15}, and others)
is the one of finding explicit formulas giving the ``general solution'' of
an under-determined system of ODEs as functions of some arbitrary functions of
time and a certain number of their time-derivatives (in fact \cite{Zerv32} allows
to change the independent variable, but we keep it to be time).
Let us call such formulas a \emph{Monge parameterization}, its \emph{order}
being the number of time-derivatives.

The authors of \cite{Flie-Lev-Mar-R92cras} already made the link with the above
mentioned work on under-determined systems of ODEs dating
back from the beginning of 20\textsuperscript{th} century; for instance, they used
\cite{Hilb12,Cart15} to obtain, in \cite{Rouc94,Mart-Rou94} some results on
flatness or linearizability of control systems.

Let us precise the relation between flatness and Monge parameterizability~:
flatness is  existence of some functions
---we call this collection of functions a \emph{flat output}---
of the state, the controls and a certain number $\ordreplat$ of
time-derivatives of the control, that ``invert'' the formulas of a Monge
parameterization, i.e. a solution $t\mapsto(x(t),u(t))$ of the control system
corresponds to only one choice of the arbitrary functions of time appearing
in the parameterization, given by these functions.
Let us call $\ordreplat$ the \emph{order} of the flat output.

Characterizing differential flatness, or dynamic
state feedback linearizability is still an open problem~\cite{Flie-Lev-Mar-R99open}, apart from the
case of single-input systems~\cite{Char-Lev-Mar89,Cart15}.
The main difficulty is that the order of a parameterization or
a flat output, if there exists any, is not known beforehand:
for a given system, if one can construct a parameterization, or a flat output,
it has a definite order, but if, for some integer $j$, one prove that there is
no parameterization of order $j$, then it might admit a parameterization of
higher order, and we do not know any a priori bound on the possible $j$'s.
In the present paper, we consider systems of the smallest dimensions for which
the answer is not known; we do not really overcome the above mentioned ``main difficulty'', in the sense that
we only say that our class of systems does not admit a parameterization of
order less than some numbers, but the description of the parameterization
that we give, and the resulting system of PDEs is valid at any order.

Consider a general control-affine system in $\RR^4$ with two controls, where $\xi\in\RR^4$ is the state, $\widetilde{w}_1$ and $\widetilde{w}_2$ are the two scalar controls
and $X_0$, $X_1$ and $X_2$ are three smooth vector fields~:
\begin{displaymath}
  \dot{\xi}=\XX_0(\xi)+\widetilde{w}_1 X_1(\xi)+\widetilde{w}_2 X_2(\xi)\ .
\end{displaymath}
In \cite{Pome97cocv}, one can find a
necessary and sufficient condition on $X_0$, $X_1$, $X_2$ for this system to admit a
flat output depending on the state and control only ($\ordreplat=0$ according to
the above notations).
Systems who do \emph{not} satisfy this conditions may or may not
admit flat outputs depending also on some time-derivatives of the control
($\ordreplat>0$). This is recalled and commented in section~\ref{sec-plan} and \ref{sec-main}.

Instead of the above control system, we study a reduced equation
(\ref{sys3}); let us briefly explain why it represents, modulo a possibly
dynamic feedback transformation, all the relevant cases.
Systems for which the iterated Lie brackets
of $X_1$ and $X_2$ do not have maximum rank can be treated in a rather simple
manner \cite[first cases of Theorem~3.1]{Pome97cocv}; if on the contrary
iterated Lie brackets do have maximum rank, it is well known (Engel normal form for
distributions of rank 2 in $\RR^4$, see~\cite{Brya-Che-Gar-G-G91}) that,
after a nonsingular feedback
($\widetilde{w}_i=\beta^{i,0}(\xi)+\beta^{i,1}(\xi)w_1+\beta^{i,2}(\xi)w_2$, $i=1,2$,
with $\beta^{1,1}\beta^{2,2}-\beta^{1,2}\beta^{2,1}\neq0$),
 there are
coordinates such that the system reads
\begin{equation}
\label{sys4}
\dot{\xi}_1=w_1\,,\ \ 
\dot{\xi}_2=\gamma(\xi_1,\xi_2,\xi_3,\xi_4)+\xi_3w_1\,,\ \ 
\dot{\xi}_3=\delta(\xi_1,\xi_2,\xi_3,\xi_4)+\xi_4w_1\,,\ \ 
\dot{\xi}_4=w_2
\end{equation}
with some smooth functions $\gamma$ and $\delta$.
One can eliminate $w_1$ and $w_2$ and, renaming $\xi_1,\xi_2,\xi_3,\xi_4$ as
$x,y,z,w$, obtain the two following relations between these four functions of time~:
\begin{equation} 
\label{sysbis}
\dot{y}=\gamma(x,y,z,w)+z\dot{x}\,,\ \ \ 
\dot{z}=\delta(x,y,z,w)+w\dot{x}
\end{equation}
(this can also be seen as a control system with state $(x,y,z)$ and controls
$w$ and $\dot{x}$).
If $\gamma$ does not depend on $w$, this system is always
parameterizable, and even flat (see \cite{Pome97cocv} or Example~\ref{ex-x} below).
If, on the contrary, $\gamma$ does depend on its last argument, one can,
around a point where the partial derivative is nonzero, invert $\gamma$ with
respect to $w$, \ie\  transform the first equation into
$w=g(x,y,z,\dot{y}-z\dot{x})$ for some function $g$, and obtain, substituting
into the last equation, a single differential relation between
$x,y,z$ written as (\ref{sys3}) in next section.

Note that (\ref{sys3}) also represents the general (non-affine) systems in
$\RR^3$ with two controls that satisfy the necessary condition given in
\cite{Rouc94,Slui93}, i.e. they are ``ruled''; we do not develop this here, see
\cite{Avan05th} or a future publication.

The paper technically focuses on Monge parameterizations of (\ref{sys3}).
The problem is unsolved if $g$ and $h$ are such that system (\ref{sys4}) does
not satisfy the above mentioned necessary and sufficient condition.
We do not give a complete solution, but our results are more general than ---and
imply---
these of \cite{Pome97cocv}.
The techniques used in the present paper, derived from the original
proof of non-parameterizability of some special systems 
in~\cite{Hilb12} (see also~\cite{Rouc92}), are much simpler and elementary that these of
\cite{Pome97cocv}: recovering the results from that paper in this
way has some interest in itself.


\section{Problem statement}
\label{sec-state}

\subsection{The systems under consideration}
\label{sec-sys}

This paper studies the solutions $t\mapsto(x(t),y(t)$, $z(t))$ of the scalar
differential equation
\begin{equation}
\label{sys3}
  \dot{z}\ \;=\;\
  h(x,y,z,\lambda)\;+\;g(x,y,z,\lambda)\,\dot{x}\ \ \ \ 
\mbox{with}\ \ \ \ 
\lambda=\dot{y}-z\dot{x}\ 
\end{equation}
where $g$ and $h$ are two real analytic functions   $\Omega\to\RR$, $\Omega$ being
an open connected subset of $\RR^4$.  
We assume that $g$ does depend on $\lambda$; more precisely, associating
to $g$ a map
  $G:\Omega\to\RR^4$ defined by 
    $G(x,y,z,\lambda)=(x,y,z,g(x,y,z,\lambda))$,
and
denoting by $g_4$ the partial derivative of $g$ with respect
  to its fourth argument,
  \begin{equation}
    \label{g4}
    g_4\mbox{ does not vanish on } \Omega
\ \ \ \mbox{and}\ \ \ 
    G\mbox{ defines a diffeomorphism }\Omega\to G(\Omega)\;.
  \end{equation}
We denote by $\widehat{\Omega}$ the open connected subset of $\RR^5$ defined from
  $\Omega$ by~:
  \begin{equation}
    \label{OO}
    (x,y,z,\dot{x},\dot{y})\in\widehat{\Omega}\ \Leftrightarrow\
    (x,y,z,\dot{y}-z\dot{x})\in\Omega\ .
  \end{equation}
From $g$ and $h$ one may define $\gamma$ and $\delta$, two real analytic
functions $G(\Omega)\to\RR$, such that
$G^{-1}(x,y,z,w)=(x,y,z,\gamma(x,y,z,w))$ and $\delta=h\circ G^{-1}$, \ie
\begin{eqnarray}\label{defgamma}
&\!\!\!\!w=g(x,y,z,\lambda)\Leftrightarrow \lambda=\gamma(x,y,z,w)\ ,
\\\label{defdelta}&
h(x,y,z,\lambda)=\delta(x,y,z,g(x,y,z,\lambda))\,,
\ \ \ 
\delta(x,y,z,w)=h(x,y,z,\gamma(x,y,z,w))\,.
\end{eqnarray}
Then, one may associate to (\ref{sys3}) the control-affine
system (\ref{sys4}) in $\RR^4$ with two controls, that can also be written as
(\ref{sysbis}); our interest however focuses on
system (\ref{sys3}) defined by $g$ and $h$ as above. Let us set some conventions~:
\begin{description}
\item[The functions $\gamma$ and $\delta$] when using the notations $\gamma$ and $\delta$, \emph{it
  is not assumed that they are related to $g$ and $h$} by
  (\ref{defgamma}) and (\ref{defdelta}), unless this is explicitly stated.

\item[Notations for the derivatives] 
  We denote partial derivatives by subscript
  indexes. For functions of
  many variables, like $\varphi(u,\ldots,u^{(k)},v,\ldots,v^{(\ell)})$ in
  (\ref{para1}), we use the name of the variable as a subscript~: 
  $p_{xu^{(k-1)}}$ means $\partial^2p/\partial x \partial u^{(k-1)}$, $\varphi_{v^{(\ell)}}$ means $\partial\varphi/\partial v^{(\ell)}$
   in
  (\ref{EDPp}-b). 
  Since the arguments of $g$, $h$, $\gamma$, $\delta$ and a few other
  functions will sometimes be 
  intricate functions of other variables, we use numeric subscripts for their
  partial derivatives~:
  $h_2$ stands for $\partial h/\partial y$, or $g_{4,4,4}$ for
  $\partial^3 g/\partial\lambda^3$. 
  To avoid confusions, we will not use numeric subscripts for other purposes
  than partial derivatives, except the subscript 0, as in
  $(x_0,y_0,z_0,\dot{x}_0,\dot{y}_0)$ for a reference point. 

  The dot denotes, as usual, derivative with respect to time, and
  $^{(j)}$ the $j$\textsuperscript{th} time-derivatives. 
\end{description}

\smallskip
The following elementary lemma 
---we do write it for the argument is used repeatedly throughout
the paper---
states that no differential equation
independent from (\ref{sys3}) can be satisfied identically by \emph{all} solutions of
(\ref{sys3})~:
\begin{lmm}
  \label{lem-absurd}
For $M\in\NN$, let $W$ be an open subset of $\RR^{3+2M}$ and
$R:W\to\RR$ a smooth function.
If \emph{any} solution $(x(.),y(.),z(.))$ of system (\ref{sys3}), 
defined on some time-interval $I$ 
and such that 
$(z(t),x(t),\ldots,x^{(M)}(t),y(t),\ldots$ $y^{(M)}(t))$ is in $W$ for all $t$
in $I$, satisfies
$\ \ R(z(t),y(t),\ldots,y^{(M)}(t),x(t),\ldots,x^{(M)}(t))=0\ $ identically on
$I$, then $R$ is identically zero on $W$.
\end{lmm}
\begin{proof}
  For \emph{any} $\mathcal{X}\in W$ there is a germ of solution of (\ref{sys3})
  such that
  $(z(0),x(0),\ldots,x^{(M)}(0)$, $y(0),\ldots,y^{(M)}(0))=\mathcal{X}$. Indeed,
  take e.g. for $x(.)$ and $y(.)$ the polynomials in $t$ of degree $M$ that
  have these derivatives at time zero; Cauchy-Lipschitz theorem then yields a
  (unique) $z(.)$ solution of (\ref{sys3}) with the prescribed $z(0)$.  
\end{proof}

\subsection{The notion of parameterization}
In order to give rigorous definitions without taking care of
time-intervals of definition of the solutions, we consider germs of solutions
at time 0, instead of solutions themselves. 
For $O$ an open subset of $\RR^n$, the notation
$\mathcal{C}_0^{\infty}({\mathbb{R}},O)$ stands for the set of germs at $t=0$
of smooth functions 
of one variable with values in $O$, see e.g. \cite{Golu-Gui73}.

Let $k$, $\ell$, $L$ be some non negative integers, $U$ an open subset of
${\mathbb{R}}^{k+\ell+2}$ and
$V$ an open subset of $\RR^{2L+3}$.
We denote by
${\mathcal{U}}\subset{\mathcal{C}}_0^{\infty}({\mathbb{R}},{\mathbb{R}}^2)$
(resp. $\mathcal{V}\subset{\mathcal{C}}_0^{\infty}({\mathbb{R}},{\mathbb{R}}^3)$~)
the set of germs
of smooth functions $t\mapsto(u(t),v(t))$ (resp. $t\mapsto(x(t),y(t),z(t))$~)
such that their jets at $t=0$ to the order precised below are
in $U$ (resp. in $V$)~:
\begin{eqnarray}\label{defu}
{\mathcal{U}}=\{(u,v)\in {\mathcal{C}}_0^{\infty}(\mathbb{R},\mathbb{R}^2) |
(u(0),\dot{u}(0),\ldots,u^{(k)}(0),v(0),\ldots,v^{(\ell)}(0)) \in U \},
\\\label{defv}
\mathcal{V}=\{(x,y,z)\in\mathcal{C}_0^{\infty}(\RR,\RR^3) |
(x(0),y(0),z(0),\dot{x}(0),\dot{y}(0),\ldots,x^{(L)}(0),y^{(L)}(0)) \in V \}.
\end{eqnarray}
These are open sets for the Whitney
${\mathcal{C}}^{\infty}$ topology \cite[p.~42]{Golu-Gui73}.

\begin{dfntn}[Monge parameterization]\label{defparam}
Let $k,\ell, L$ be non negative integers, $L>0$, $k\leq\ell$, and
$\mathcal{X}=(x_0,y_0,z_0,\dot{x}_0,\dot{y}_0,\ldots,x^{(L)}_0$, $y^{(L)}_0)$ a point
in $\widehat{\Omega}\times\RR^{2L-2}$ ($\widehat{\Omega}$ is defined in (\ref{OO})).
A \emph{parameterization of order $(k,\ell)$ at $\mathcal{X}$} for system
(\ref{sys3}) is defined by
\begin{itemize}
\item 
a neighborhood $V$ of $\mathcal{X}$ in $\widehat{\Omega}\times\RR^{2L-2}$,
\item 
an open subset $U\subset{\mathbb{R}}^{k+\ell+2}$ and
\item 
three real analytic functions $U\to {\mathbb{R}}$, denoted
$\varphi$, $\psi$, $\chi$,
\end{itemize}
such that, with $\mathcal{U}$ and $\mathcal{V}$  defined
from $U$ and $V$ according
to (\ref{defu})-(\ref{defv}), and\\
$\Gamma:\mathcal{U}\to{\mathcal{C}}_0^{\infty}(\mathbb{R},\mathbb{R}^3)$ the
map that assigns to $(u,v)\in{\mathcal{U}}$ the germ $\Gamma(u,v)$ at $t=0$ of
\begin{equation}
\label{para1}
t\ \mapsto\ \left(\begin{array}{c}x(t)\\y(t)\\z(t)\end{array}\right)
=
\left(
\begin{array}{l}
\varphi(u(t),\dot{u}(t),\ldots,u^{(k)}(t),v(t),\dot{v}(t),\ldots,v^{(\ell)}(t))\\
\psi(u(t),\dot{u}(t),\ldots,u^{(k)}(t),v(t),\dot{v}(t),\ldots,v^{(\ell)}(t))\\
\chi(u(t),\dot{u}(t),\ldots,u^{(k)}(t),v(t),\dot{v}(t),\ldots,v^{(\ell)}(t))
\end{array}\right)
\ ,
\end{equation}
the following three properties hold~:
\begin{enumerate}
\item\label{def2}
 for all $(u,v)$ belonging to $\mathcal{U}$,
$\Gamma(u,v)$
is a solution of system (\ref{sys3}),
\item\label{def1}
the map $\Gamma$ is open and $\Gamma({\mathcal{U}})\supset{\mathcal{V}}$,
\item\label{def3}
the two maps $U\to{\mathbb{R}}^3$ defined by the triples $(\varphi_{u^{(k)}},\psi_{u^{(k)}},\chi_{u^{(k)}})$ and\\
$(\varphi_{v^{({\ell})}},\psi_{v^{({\ell})}},\chi_{v^{({\ell})}})$ are
identically zero on no open subset of $U$.
\end{enumerate}
\end{dfntn}
\begin{rmrk}[On ordering the pairs $(k,\ell)$]
  \label{rmk-kl}
Since $u$ and $v$ play a symmetric role, they can always be exchanged, and
there is no lack of generality in assuming $k\leq\ell$. 
This convention is useful only when giving bounds on $(k,\ell)$. For
instance, $k\geq2$ means that both integers are no smaller than 2.
\end{rmrk}

\begin{xmpl}
\label{ex-xu}
  Consider the equation $\ \dot{z}=y+(\dot{y}-z\dot{x})\dot{x}\ $, \ie\
  (\ref{sys3}) with $g=\lambda$, $h=y$ (and $\widehat{\Omega}=\RR^3$). At any
  $(x_0,y_0,z_0,\dot{x}_0,\dot{y}_0,\ddot{x}_0,\ddot{y}_0)$ such that
  $\ddot{x}_0+{\dot{x}_0}^{\,3}\neq1$, a parameterization of order $(1,2)$ is
  given by~:
  \begin{equation}
    \label{eq:ex-xu}
    x=v\,,\ \ 
    y=\frac{{\dot{v}}^2u+\dot{u}}{\ddot{v}+{\dot{v}}^3-1}\,,\ \ 
    z=\frac{(1-\ddot{v})u+\dot{v}\dot{u}}{\ddot{v}+{\dot{v}}^3-1}\,.
  \end{equation}
It is easy to check that $(x,y,z)$ given by these formulas does satisfy the
equation, point~\ref{def1} is true because the above formulas can be
``inverted'' by $u=-z+y\dot{x}$, $v=x$ (this gives the ``flat output'' see
section~\ref{sec-flat}), point~\ref{def3} is true because $\psi_{\dot{u}}$,
$\psi_{\ddot{v}}$, $\chi_{\dot{u}}$ and $\chi_{\ddot{v}}$ are nonzero rational
functions. Here, $L=2$ and $V$ can be taken the whole set of
$(x,y,z,\dot{x},\dot{y},\ddot{x},\ddot{y})\in\RR^7$ such that
$\ddot{x}+{\dot{x}}^{3}\neq1$ and $U$ the whole set of
$(u,\dot{u},v,\dot{v},\ddot{v})\in\RR^5$ such that
$\ddot{v}+{\dot{v}}^3\neq1$. 
\end{xmpl}

\begin{xmpl}
\label{ex-x}
  Suppose that the function $\gamma$ in (\ref{sysbis}) depends on $x,y,z$ only
  (this is treated in \cite[case 6 in Theorem~3.1]{Pome97cocv}).
  For such systems, eliminating $w$ does not lead to (\ref{sys3}),
  but to the simpler relation $\dot{y}-z\dot{x}=\gamma(x,y,z)$. One can easily
  adapt the above definition replacing (\ref{sys3}) by this relation.
  This system $\dot{y}-z\dot{x}=\gamma(x,y,z)$ admits a
  parameterization of order (1,1) at any $(x_0,y_0,z_0,\dot{x}_0,\dot{y}_0)$ such
  that $\dot{x}_0+\gamma_3(x_0,y_0,z_0)\neq 0$.
  \\\textit{Proof.}\ 
      In a neighborhood of such a point,, the map
  $(x,\dot{x},y,z)\mapsto(x,\dot{x},y,\gamma(x,y,z)+z\dot{x})$ is a 
  local diffeomorphism, whose inverse can be written as
  $(x,\dot{x},y,\dot{y})\mapsto(x,\dot{x},y$, $\chi(x,\dot{x},y,\dot{y}))$, thus
  defining a map $\chi$.
  Then $x=u, y=v, z=\chi(u,\dot{u},v,\dot{v})$ defines a parameterization of order
  (1,1) in a neighborhood of these points. 
\end{xmpl}

\begin{rmrk}
\label{rmrk-L}
The integer $L$ characterizes the number of derivatives needed
to describe the open set where the parameterization is valid.
For instance, in Examples~\ref{ex-xu} and \ref{ex-x}, $L$ must be taken no smaller than
2 and 1 respectively. 
Obviously, a parameterization of order $(k,\ell)$ at
$(x_0,y_0,z_0,\dot{x}_0,\dot{y}_0,\ldots,x^{(L)}_0,y^{(L)}_0)$ is also, for $L'>L$ and 
\emph{any} $(x^{(L+1)}_0,y^{(L+1)}_0,\ldots,x^{(L')}_0,y^{(L')}_0)$,
a parameterization of the same
order at $(x_0,y_0,z_0,\dot{x}_0,\dot{y}_0,\ldots,x^{(L')}_0,y^{(L')}_0)$.
\end{rmrk}

The above definition is local around some jet of solutions of
(\ref{sys3}). In general, the idea of a global parameterization, meaning
that $\Gamma$ would be defined globally, is not realistic; it is not 
realistic either to require that there exists a parameterization around all jets
(this would be ``everywhere local'' rather than ``global'')~: the systems in
example~\ref{ex-x} admit a local parameterization around ``almost every'' jets,
meaning jets outside the zeroes
of a real analytic function (namely jets such that $\dot{x}+\gamma_3(x,y,z)\neq0$).
We shall not define more precisely the notion of ``almost everywhere local''
parameterizability, but rather the following (sloppier) one.
\begin{dfntn}
\label{def-param-bof}
We say that system (\ref{sys3}) \emph{admits a parameterization of
  order $(k,\ell)$ somewhere in $\Omega$} if there exist an integer $L$ and at
least one jet $(x_0,y_0,z_0,\dot{x}_0,\dot{y}_0,\ldots$, $x^{(L)}_0,y^{(L)}_0)\in \widehat{\Omega}\times\RR^{2L-2}$
with a parameterization of order $(k,\ell)$ at this jet in the sense of Definition~\ref{defparam}
\end{dfntn}
In a colloquial way this is a ``somewhere local'' property. Using real
analyticity, it should imply ``almost everywhere local'', but
we do not investigate this.

\subsection{The functions $S$, $T$ and $J$}
\label{sec-ST}
Given $g,h$,
let us define three functions $S$, $T$ and $J$, to be used to discriminate
different cases. They were already more or less present in \cite{Pome97cocv}.
The most compact way is as follows~: let $\omega$, $\omega^1$ and $\eta$ be the following
differential forms in the variables $x,y,z,\lambda$~: 
\begin{equation}
\label{om}
\begin{array}{ll}
\!\!\!\!\omega^1=\xdif{y}-z\xdif{x}\,,
&
\omega \ =\  -2\,{g_{4}}^{2}\xdif{x}
+\left(g_{4,4}\,h_{4}- g_{4}\,h_{4,4}\right)\omega^1
-g_{4,4}\left(\xdif{z}-g\xdif{x}\right),
\\
& 
\eta \ = \ \xdif{z}-g\xdif{x}-h_4\,\omega^1\ .
\end{array}
\end{equation}
From (\ref{g4}), 
$\omega\!\wedge\!\omega^1\!\wedge\!\eta=2{g_4}^2\xdif{x}\!\wedge\!\xdif{y}\!\wedge\!\xdif{z}\neq0$.
Decompose $\xdif{\omega}\!\wedge\!\omega$ on the basis
$\omega,\omega^1,\eta,\xdif{\lambda}$, thus defining the functions $S$, $T$
and $J$ (we say more on their expression and meaning in section~\ref{sec-ST0})~:
\begin{eqnarray}
  \label{STJ}
  \xdif{\omega}\wedge\omega &=& -\left(
\frac{S}{2g_4}\,\xdif{\lambda}\wedge\eta+\frac{T}{2}\,\xdif{\lambda}\wedge\omega^1
\;+\;J\,\omega^1\wedge\eta
\right)\wedge\omega\ .
\end{eqnarray}

\begin{xmpl}\label{3ex}
Le us illustrate the computation of $S$, $T$
and $J$ on the following three particular cases of (\ref{sys3}). For each of
them, the table below gives the differential forms $\omega$ and $\eta$, the
decomposition of $\xdif{\omega}\wedge\omega$ on $\omega^1, \omega, \eta,
\xdif{\lambda}$ and the resulting $S,T,J$ according to (\ref{STJ}).
System (a) was already studied in Example~\ref{ex-xu}.
\begin{equation}
  \label{eq:ex}
  \mbox{(a):}\ \dot{z}=y+(\dot{y}-z\dot{x})\dot{x}\;,
\ \ \ \ \ \ 
  \mbox{(b):}\ \dot{z}=y+(\dot{y}-z\dot{x})(\dot{y}-(z-1)\dot{x})\;,
\ \ \ \ \ \ 
  \mbox{(c):}\ \dot{z}=y+(\dot{y}-z\dot{x})^2\dot{x}\;.
\end{equation}
\begin{center}
\begin{tabular}{c|cc|c|c|c|}
\begin{tabular}{c} system\\(\ref{eq:ex})\end{tabular} & $g (x,y,z,\lambda)$ & $h(x,y,z,\lambda)$ & 
$ \begin{array}{c} -\omega/2 \\ \eta \end{array}$ 
& $\xdif{\omega}\wedge\omega$
& $S,T,J$
\\
\hline
(a)
& $\lambda$ & $y$
& $\begin{array}{l}
    \xdif{x}
    \\
    \xdif{z}-\lambda\xdif{x}
    \end{array}   $
& 0 & 0, 0, 0
\\
\hline
(b)
& $\lambda$ & $y+\lambda^2$
& $\begin{array}{l}
    \xdif{y}-(z-1)\xdif{x}
    \\
    \xdif{z}-\lambda\xdif{x}-2\lambda\omega^1
    \end{array}   $
& $\omega^1\wedge\eta\wedge\omega$
& 0, 0, $-1$
\\
\hline
(c)
& $\lambda^2$ & $y$
& $\begin{array}{l}
    \xdif{z}+3\lambda^2\xdif{x}
    \\
    \xdif{z}-\lambda^2\xdif{x}
    \end{array}   $
& $\frac{3}{\lambda}\xdif{\lambda}\wedge\eta\wedge\omega$
& $-12$, 0, 0
\\
\hline
\end{tabular}
  \end{center}
\end{xmpl}

\subsection{Contributions and organization of the paper}
\label{sec-plan}

If $S=T=J=0$, \ie\ $\xdif{\omega}\wedge\omega=0$, system (\ref{sys3}) admits a
parameterization of order (1,2), at all points except some singularities.
This is stated further as Theorem~\ref{prop-plats}, but was already contained
in~\cite{Pome97cocv}. 
We conjecture that these systems are the only parameterizable ones of
these dimensions, \ie\ system (\ref{sys3}) admits no parameterization of any
order if $(S,T,J)\neq(0,0,0)$, \ie\ if $\xdif{\omega}\wedge\omega\neq0$.

This is unfortunately still a conjecture, but we give the following results, valid if
$(S,T,J)\neq(0,0,0)$
(recall that $k\leq\ell$, see Remark~\ref{rmk-kl})~:
\begin{itemize}
\item 
system (\ref{sys3}) admits no parameterization of order $(k,\ell)$ with $k\leq
2$ or $k=\ell=3$ (Theorem~\ref{th-3}),
\item 
a parameterization of order
$(k,\ell)$ must come from a solution of the system of PDEs $\sysedp$ (Theorem~\ref{th-cns}),
\item since a solution of this system of PDEs is also sufficient to construct a
parameterization (Theorem~\ref{edpparam}), the conjecture can be entirely re-formulated in terms
of this system of partial differential relations.
\end{itemize}
Note that this allows one to recover the results
from \cite{Pome97cocv} on $(x,u)$-flatness\footnote{The term ``dynamic
  linearizable'' in \cite{Pome97cocv} is synonymous to ``flat'' here.
}. See Remark~\ref{rmrk-old} for details.

\smallskip
The paper is organized as follows.
Section~\ref{sec-edp} is about the above mentioned partial differential system $\sysedp$.
Section~\ref{sec-ST0} is devoted to some special constructions for the case
where $S=T=0$, and geometric interpretations.
The main results are stated in Section~\ref{sec-main}, based on sufficient
conditions obtained in Sections~\ref{sec-edp} and \ref{sec-ST0},  and necessary
conditions stated and proved in Section~\ref{sec-nec}.
Sections~\ref{sec-flat} and \ref{sec-concl} comment on flatness vs. Monge
parameterization and then give a conclusion and perspectives.

\section{A system of partial differential equations}
\label{sec-edp}

\emph{This section can profitably be skipped or overlooked in a first
  reading}; the reader will come back when needed to this material that might
appear, at first sight, somehow disconnected from the thread of the paper. 

It defines $\sysedp$ and its ``regular solutions'', proves that a regular solution induces a
parameterization of order $(k,\ell)$, and that no regular solution exists
unless $k\geq3$ and $\ell\geq4$.

\subsection{The equation $\sysedp$, regular solutions}

For $k$ and $\ell$ some positive integers, we define a partial differential
system in  $k+\ell+1$ independent variables and one dependent variable,
\ie\  the unknown is one function of $k+\ell+1$ variables. 
The dependent variable is denoted by $p$ and the independent variables by
$u,\dot{u},\ldots,u^{(k-1)},x,v,\dot{v},\ldots,v^{({\ell-1})}$.
Although the names of the variables may suggest
``time-derivatives'', time is \emph{not} a variable here.

In $\RR^{k+\ell+1}$ with the independent variables as coordinates, let $F$ be the differential operator of order 1
\begin{equation}
\label{F}
F=\sum_{i=0}^{k-2}u^{(i+1)}
\partialx{u^{(i)}}+\displaystyle\sum_{i=0}^{\ell-2}v^{(i+1)}\partialx{v^{(i)}}\
,
\end{equation}
where the first sum is zero if $k\leq 1$ and the second one is zero if $\ell\leq 1$.

Let $\widetilde{\Omega}$ be an open connected subset of $\RR^4$ and 
$\gamma,\delta$ two real analytic functions $\widetilde{\Omega}\to\RR$ such that
$\gamma_4$ (partial derivative of $\gamma$ with respect to its
4\textsuperscript{th} argument, see end
of section~\ref{sec-sys}) does not vanish on
$\widetilde{\Omega}$.
Consider the system of  two partial differential equations and three inequations~:
\begin{equation}
\label{EDPp}
\sysedp\left\{
\begin{array}{ll}
p_{u^{(k-1)}}\bigl(Fp_x-\delta(x,p,p_x,p_{xx})\bigr)
-p_{xu^{(k-1)}}\bigl(Fp-\gamma(x,p,p_x,p_{xx})\bigr)=0\;,&\mbox{(a)}\\
p_{u^{(k-1)}}\,p_{xv^{(\ell-1)}}-p_{xu^{(k-1)}}\,p_{v^{(\ell-1)}}=0\;,
&\mbox{(b)}
\\
p_{u^{(k-1)}}\neq 0\;,&\mbox{(c)}
\\
p_{v^{(\ell-1)}}\neq 0\;,&\mbox{(d)}\\
\gamma_1+\gamma_2\,p_x+\gamma_3\,p_{xx}+\gamma_4\,p_{xxx}-\delta\neq0\;.
&\mbox{(e)}
\end{array}\right.\!\!\!\end{equation}
To any $p$ satisfying $\sysedp$, we associate
two functions $\sigma$ and $\tau$, and a vector field $E$~:
\begin{equation}\label{defsigtau}
\sigma=-\frac{p_{v^{(\ell-1)}}}{p_{u^{(k-1)}}}
\;,\ 
\tau=\frac{-Fp+\gamma(x,p,p_x,p_{xx})}{p_{u^{(k-1)}}}\;,\ E=\sigma\partialx{u^{(k-1)}}+\partialx{v^{(\ell-1)}}\,.
\end{equation}
We also introduce the differential operator $D$ (see Remark~\ref{rmk-spur} on
the additional variables $\dot{x},\ldots,x^{(k+\ell-1)}$)~:
\begin{equation}
\label{D}
D=F+\tau\partialx{u^{(k-1)}}+
\displaystyle\sum_{i=0}^{k+\ell-2}x^{(i+1)}\partialx{x^{(i)}}\ .
\end{equation}
\begin{dfntn}[Regular solutions of $\sysedp$]
\label{def-Kreg}
A \emph{regular} solution of system $\sysedp$ is a real analytic
function $p:\ouvedp\to\RR$, with $\ouvedp$ a connected open subset of
$\RR^{k+\ell+1}$, 
such that the image of $\ouvedp$ by $(x,p,p_x,p_{xx})$ is contained in
$\widetilde{\Omega}$, (\ref{EDPp}-a,b) are identically satisfied on $\ouvedp$, the left-hand
sides of (\ref{EDPp}-c,d,e) are not identically zero, 
and, for at least one integer $K\in\{1,\ldots,k+\ell-2\}$,
\begin{equation}
  \label{EDPpK}
  ED^{K}p\neq0
\end{equation}
(not identically zero, as a function of
$u,\ldots,u^{(k-1)},x,v,\ldots,v^{(\ell-1)},\dot{x},\ldots,x^{(K)}$ on
$\ouvedp\times\RR^K$). 
\hspace{2em}
We call it \emph{$K$-regular} if $K$ is the smallest such
integer, \ie\  if $ED^ip=0$ for all $i\leq K-1$.
\end{dfntn}
\begin{rmrk}[on the additionnal variables $\dot{x},\ldots,x^{(k+\ell-1)}$ in $D$]
  \label{rmk-spur}
These variables appear in the expression (\ref{D}). Note that $D$ is only
applied (recursively) to functions of
$u,\ldots,u^{(k-1)},x$, $v,\ldots,v^{(\ell-1)}$ only; hence we view it as a
vector field in $\RR^{k+\ell+1}$ with these variables as parameters.
In fact, $D$ is only used in $ED^ip$, $1\leq i\leq k+\ell-1$. This is a
polynomial with respect to the variables $\dot{x},\ddot{x},\ldots,x^{(i)}$ with coefficients depending on
$u,\ldots,u^{(k-1)},x,v,\ldots,v^{(\ell-1)}$ via the functions $p$, $\gamma$,
$\delta$ and their partial derivatives. Hence $ED^ip=0$ means that all these
coefficients are zero, \ie\ it encodes a collection of differential
relations on $p$, where the spurious variables
$\dot{x},\ddot{x},\ldots,x^{(i)}$ no longer appear. Likewise, $ED^ip\neq 0$ means that one of these relations is
not satisfied.
\end{rmrk}
\begin{dfntn}
\label{def-Kreg-gen}
We say that system $\sysedp$ \emph{admits a regular (resp. $K$-regular) solution somewhere in
  $\widehat{\Omega}$}
if there exist at least an open connected $\ouvedp\subset\RR^{k+\ell+1}$ and a regular (resp. $K$-regular) solution $p:\ouvedp\to\RR$.
\end{dfntn}


\begin{rmrk}
\label{rmk-EDPbis}
It is easily seen that $p$ is solution of $\sysedp$ if and only if
there exist $\sigma$ and $\tau$ such that $(p,\sigma,\tau)$ is a solution of
\begin{equation}
  \label{EDPbis}
  \begin{array}{rll}
Fp+\tau p_{u^{(k-1)}}=\gamma(x,p,p_x,p_{xx})&\hspace{1em}&
Ep=0\,,\ \sigma_x=0\,,
\\
Fp_x+\tau p_{x,u^{(k-1)}}=\delta(x,p,p_x,p_{xx})&&p_{u^{(k-1)}}\neq0\,,\
\tau_x\neq0\,,\ \sigma\neq0 
  \end{array}
\end{equation}
Indeed, (\ref{EDPp}) does imply the above relations with $\sigma$ and $\tau$
given by (\ref{defsigtau}); in particular,
$\tau_x\neq0$ is equivalent to (e) and $\sigma\neq0$ to (d); conversely,
eliminating $\sigma$ and $\tau$ in (\ref{EDPbis}), one recovers $\sysedp$.
Note also that, with $g$ and $h$ related to $\gamma$ and $\delta$ by (\ref{defgamma}) and
(\ref{defdelta}), any solution of the above equations and inequations satisfies
\begin{equation}
  \label{psolgh}
  Dp_x=h(x,p,p_x,Dp-p_x\dot{x})+g(x,p,p_x,Dp-p_x\dot{x})\dot{x}\ .
\end{equation}
\end{rmrk}

The following will be used repeatedly in the paper~:
\begin{lmm}
\label{lemAAA}
If $p$ is a solution of system $\sysedp$ and
\begin{enumerate}
\item\label{pt1} either it satisfies a relation of the type $p_x=\alpha(x,p)$ with
  $\alpha$ a function of two variables,
\item\label{pt2} or it satisfies a relation of the type $p_{xx}=\alpha(x,p,p_x)$ with
  $\alpha$ a function of three variables,
\item\label{pt3} or it satisfies two relations of the type
  $p_{xxx}=\alpha(x,p,p_x,p_{xx})$ and\\
  $Fp_{xx}+\tau p_{xxu^{(k-1)}}=\psi(x,p,p_x,p_{xx})$, with $\psi$
  and $\alpha$ two functions of four variables, 
\end{enumerate}
then it satisfies $ED^ip=0$ for all $i\geq0$ and
hence is \emph{not} a regular solution of $\sysedp$.
\end{lmm}
\begin{proof}
Point~\ref{pt1} implies point~\ref{pt2} because
differentiating the relation $p_x=\alpha(x,p)$ with respect to $x$
yields $p_{xx}=\alpha_x(x,p)+p_x\alpha_p(x,p)$.
Likewise, point~\ref{pt2} implies point~\ref{pt3}~:
differentiating the relation $p_{x,x}=\alpha(x,p,p_x)$ with respect to $x$
yields
$p_{xxx}=\alpha_x(x,p,p_x)+p_x\alpha_p(x,p,p_x)+p_{xx}\alpha_{p_x}(x,p,p_x)$
while differentiating it along the vector field 
$F+\tau\,\partial\!/\partial u^{(k-1)}$ and using (\ref{EDPbis}) yields \\
$Fp_{xx}+\tau p_{xxu^{(k-1)}}=\gamma(x,p,p_x,p_{xx})\,\alpha_p(x,p,p_x) +
\delta(x,p,p_x,p_{xx})\,\alpha_{p_x}(x,p,p_x)$. 

Let us prove that point~\ref{pt3} implies $ED^ip=ED^ip_x=ED^ip_{xx}=0$ for all
$i\geq0$, hence the lemma. It is indeed true for $i=0$ and the following three
relations:
$Dp=\gamma(x,p,p_x,p_{xx})+\dot{x}\,p_{x}$,
$Dp_x=\delta(x,p,p_x,p_{xx})+\dot{x}\,p_{xx}$,
$Dp_{xx}=\psi(x,p,p_x,p_{xx})+\dot{x}\,\alpha(x,p,p_x,p_{xx})$,
that are implied by (\ref{D}), (\ref{EDPbis}) and the two relations in point~\ref{pt3}
allow one to go from $i$ to $i+1$
($ED^ix=Ex^{(i)}=0$ and $ED^i\dot{x}=Ex^{(i+1)}=0$ from the very definition of
$D$ and $E$).
\end{proof}

\subsection{The relation with Monge parameterizations}
\label{sec-edpsuff}

Let us now explain how a Monge parameterization for system (\ref{sys3}) can be
deduced from a regular solution $p:\ouvedp\to\RR$ of $\sysedp$.
This may seem anecdotic but it is not, for we shall prove (\cf\ 
sections~\ref{sec-main} and \ref{sec-nec}) that all Monge parameterizations
are of this type, except when $g$ and $h$ are such that
$\xdif{\omega}\wedge\omega=0$ (see (\ref{om})-(\ref{STJ})).

We saw in Remark~\ref{rmk-EDPbis} that (\ref{EDPp}-e) is equivalent to
$\tau_x\neq0$; let
$(u_0,\ldots,u_0^{(k-1)},x_0,v_0,\ldots$, $v_0^{(\ell-1)})\in\ouvedp$ be such that 
$\tau_x(u_0,\ldots,u_0^{(k-1)},x_0,v_0,\ldots,v_0^{(\ell-1)})\neq0$.
Choose any $(u_0^{(k)},v_0^{(\ell)})\in\RR^2$ (for instance with
$v_0^{(\ell)}=0$) such that
\begin{equation}
  \label{eq:ukvl}
  u_0^{(k)}-\sigma(u_0,\ldots, u_0^{(k-1)},v_0,\ldots, v_0^{(\ell-1)})\,v_0^{(\ell)}
=\tau(u_0,\ldots,u_0^{(k-1)},x_0,v_0,\ldots,v_0^{(\ell-1)}).
\end{equation}
Then, the implicit function theorem provides a neighborhood
$V$ of $(u_0,\ldots,u_0^{(k)}$, $v_0,\ldots,v_0^{(\ell)})$ in $\RR^{k+\ell+2}$
and a real analytic map $\varphi:V\to\RR$ 
such that $\varphi(u_0,\ldots,u_0^{(k)}$, $v_0,\ldots,v_0^{(\ell)})=x_0$ and
\begin{equation}
\label{defphi}
\tau(u,\ldots,u^{(k-1)},\varphi(u\cdots v^{(\ell)}),v,\ldots,v^{(\ell-1)})\,=\,
u^{(k)}-\sigma(u,\ldots,u^{(k-1)},v,\ldots,v^{(\ell-1)})\,v^{(\ell)}
\end{equation}
identically on $V$. Two other maps $V\to\RR$ may be defined by
\begin{eqnarray}
\label{defpsi}
\psi(u,\ldots,u^{(k)},v,\ldots,v^{(\ell)})&=&p(u,\ldots,u^{(k-1)},\varphi(\cdots),v,\ldots,v^{(\ell-1)}),
\\
\label{defchi}
\chi(u,\ldots,u^{(k)},v,\ldots,v^{(\ell)})&=&
p_x(u,\ldots,u^{(k-1)},\varphi(\cdots),v,\ldots,v^{(\ell-1)}).
\end{eqnarray}

From these $\varphi$, $\psi$ and $\chi$, one can define a map $\Gamma$ as in
(\ref{para1}) that is a candidate for a parameterization. 
We prove below that, if $p$ is a regular solution of $\sysedp$, then this
$\Gamma$ is
indeed a parameterization, at least away from some singularities. The following lemma
describes these singularities;
it is proved in Appendix~\ref{app-1}.
\begin{lmm}\label{jacob}
Let $\ouvedp$ be an open connected subset of $\RR^{k+\ell+1}$ and $p:\ouvedp\to\RR$ be a
$K$-regular solution of system $\sysedp$, see (\ref{EDPp}). Define the map
$\pi:\ouvedp\times\RR^K\to\RR^{K+2}$ by
\begin{equation}
  \label{pi}
  \!\!\pi(u\cdots u^{(k-1)},x,v\cdots v^{(\ell-1)},\dot{x}\cdots x^{(K)})=\!\left(\!\!\!
  \begin{array}{c}
p_x(u\cdots u^{(k-1)},x,v\cdots v^{(\ell-1)})\\
p(u\cdots u^{(k-1)},x,v\cdots v^{(\ell-1)})
\\Dp(u\cdots u^{(k-1)},x,v\cdots v^{(\ell-1)},\dot{x})\\\vdots
\\D^Kp(u\cdots u^{(k-1)},x,v\cdots v^{(\ell-1)},\dot{x}\cdots x^{(K)})
  \end{array}\!\!\!\right).
\end{equation}
There exist two non-negative integers  $i_0\leq k$ and
$j_0\leq\ell$ such that $i_0+j_0=K+2$ and
\begin{equation}
  \label{detpi}
  \det\left(
\partialxy{\pi}{u^{(k-i_0)}},\ldots,\partialxy{\pi}{u^{(k-1)}},
\partialxy{\pi}{v^{(\ell-j_0)}},\ldots,\partialxy{\pi}{v^{(\ell-1)}}\right)
\end{equation}
is a nonzero real analytic function on $\ouvedp\times\RR^K$.
\end{lmm}

We can now state precisely the announced sufficient condition.
Its interest is discussed in Remark~\ref{rmk-suff}.

\begin{thrm}\label{edpparam}
Let $p:\ouvedp\to\RR$, with $\ouvedp\subset\RR^{k+\ell+1}$ open,  be a
$K$-regular solution of system $\sysedp$, and $i_0,j_0$ be given by Lemma~\ref{jacob}.
Then, the maps $\varphi,\psi,\chi$ constructed above define a parameterization
$\Gamma$ of system (\ref{sys3}) of order $(k,\ell)$ (see
Definition~\ref{defparam}) at any 
jet of solutions $(x_0,y_0,z_0,\dot{x}_0,\ldots,x_0^{(K)},\dot{y}_0,\ldots,y_0^{(K)})$
such that, 
for some $u_0,\ldots,u_0^{(k-1)},v_0,\ldots,v_0^{(\ell-1)}$, 
\begin{equation}\label{lien}\left.\begin{array}{l}
(u_0,\ldots,u_0^{(k-1)},x_0,v_0,\ldots,v_0^{(\ell-1)})\in\ouvedp\,,\\
z_0=p_x(u_0,\ldots, u_0^{(k-1)}, v_0,\ldots, v_0^{(\ell-1)}, x_0)\,,\\
y_0^{(i)}=D^ip(u_0,\ldots, u_0^{(k-1)}, v_0,\ldots, v_0^{(\ell-1)}, x_0, \ldots,
  x_0^{(i)})\ \ \ 0\leq i\leq K\,,
\end{array}\right\}
\end{equation}
the left-hand sides of (\ref{EDPp}-c,d,e) are all nonzero at 
$(u_0,\ldots,u_0^{(k-1)},x_0,v_0,\ldots,v_0^{(\ell-1)})$,
and the function $ED^Kp$ and the determinant (\ref{detpi}) are nonzero at point
$(u_0,\ldots,u_0^{(k-1)},x_0,\ldots,x_0^{(K)}$, $v_0,\ldots,v_0^{(\ell-1)})\in
\ouvedp\times\RR^{K}\ $.
\end{thrm}

\begin{proof}
Let us prove that $\Gamma$ given by (\ref{para1}),
with the maps $\varphi,\psi,\chi$ constructed above,
satisfies the three points of Definition~\ref{defparam}.
Differentiating (\ref{defphi}) with respect to $u^{(k)}$ and $v^{(\ell)}$ yields
$\varphi_{u^{(k)}}\tau_x=1$, $\varphi_{v^{(\ell)}}\tau_x=-\sigma$, hence the
point~\ref{def3} ($\sigma\neq0$ from (\ref{EDPbis})).
To prove point \ref{def2},
let $u(.),v(.)$ be arbitrary and $x(.),y(.),z(.)$ be defined by (\ref{para1}).
Differentiating (\ref{para1}) with respect to time, using relations
(\ref{defpsi}) and (\ref{defchi}), taking $u^{(k)}(t)$ from
(\ref{defphi}), one has
$$
\dot{y}(t)=Fp+\tau p_{u^{(k-1)}}+v^{(\ell)}(t)Ep+\dot{x}(t)\,z(t)\,,\ \ \ 
\dot{z}(t)=Fp_x+\tau p_{x,u^{(k-1)}}+v^{(\ell)}(t)Ep_x+\dot{x}(t)\,p_{xx}\,,
$$
where $F$ is given by (\ref{F}) and the argument $(u(t)\ldots u^{(k-1)}(t),x(t),v(t)\ldots
v^{(\ell-1)}(t))$ for $Fp$, $Fp_x$ $Ep$, $Ep_x$, $\tau$, $p_{x,u^{(k-1)}}$, $p_{u^{(k-1)}}$ and $p_{xx}$ is omitted.
Then, (\ref{EDPbis}) implies, again omitting
the arguments of $p_{xx}$, one has
$\dot{y}(t)=\gamma(x(t),y(t),z(t),p_{xx})+z(t)\dot{x}(t)$, and 
$\dot{z}(t)=\delta(x(t),y(t),z(t),p_{xx})+p_{xx}\dot{x}(t)$.
The first equation yields $p_{xx}=g(x(t),y(t)$, $z(t),\dot{y}(t)-z(t)\dot{x}(t))$
with $g$ related to $\gamma$ by (\ref{defgamma}), and then the second one
yields (\ref{sys3}), with $h$ related to $\delta$ by (\ref{defdelta}).
This proves point~\ref{def2}. The rest of the proof is devoted to point \ref{def1}.

Let $t\mapsto(x(t),y(t),z(t))$ be a solution of (\ref{sys3}).
We may consider $\Gamma(u,v)=(x,y,z)$ (see (\ref{para1}))
as a system of three
ordinary differential equations in two unknown functions $u, v$~:
\begin{eqnarray}
  \label{eqq101}
  \!\!\!\!u^{(k)} -\sigma(u,\ldots,u^{(k-1)},v,\ldots, v^{(\ell-1)}) v^{(\ell)} -
  \tau(u,\ldots,u^{(k-1)},x,v,\ldots, v^{(\ell-1)})
&=&0,
\\
\label{eqq102}
p(u,\ldots,u^{(k-1)},x,v,\ldots, v^{(\ell-1)}) &=&y,
\\
\label{eqq103}
p_x(u,\ldots,u^{(k-1)},x,v,\ldots, v^{(\ell-1)}) &=&z.
\end{eqnarray}
Differentiating (\ref{eqq102}) $K+1$ times, substituting $u^{(k)}$ from (\ref{eqq101}), and using the fact
that $ED^ip=0$ for $i\leq K$ (see Definition~\ref{def-Kreg}), we get
\begin{eqnarray}
\label{eqq102i}
  D^ip\,(u(t),\ldots,u^{(k-1)}(t),v(t),\ldots,
v^{(\ell-1)}(t),x(t),\ldots,x^{(i)}(t))=\frac {d^iy}{dt^i}(t)
\;,\ \ 1\leq i\leq K,
\\[1ex]
\nonumber
v^{(\ell)}(t)\;ED^{K}p\,(u(t),\ldots,u^{(k-1)}(t),v(t),\ldots,
v^{(\ell-1)}(t),x(t),\ldots,x^{(K)}(t))
\hspace{8em}
\\
\label{eqq102K}
\;+\;D^{K+1}p\,(u(t),\ldots,u^{(k-1)}(t),v(t),\ldots,
v^{(\ell-1)}(t),x(t),\ldots,x^{(K+1)}(t))  
  =\frac {d^{K+1}y}{dt^{K+1}}(t)\,.
\end{eqnarray}
Equations (\ref{eqq102})-(\ref{eqq103})-(\ref{eqq102i}) can be written
\vspace{-1.2\baselineskip}
\begin{equation}
  \label{eq:pipi}
  \pi(u,\ldots,\ldots,u^{(k-1)},x,v,\ldots,
v^{(\ell-1)},\dot{x},\ldots,x^{(K)})\ =\ 
\left(
  \begin{array}{c}
z\\y\\\dot{y}\\\vdots\\y^{(K)}
  \end{array}\right)
\end{equation}
with $\pi$ given by (\ref{pi}).
From the implicit function theorem, since the determinant (\ref{detpi}) is nonzero,
(\ref{eqq102})-(\ref{eqq103})-(\ref{eqq102i})
yields
$u^{(k-i_0)},\ldots,u^{(k-1)}$, $v^{(\ell-j_0)},\ldots,v^{(\ell-1)}$ as explicit
functions of $u,\ldots,u^{(k-i_0-1)}$, $v,\ldots,v^{(\ell-j_0-1)}$,
$x,\ldots,x^{(K)}$, $y,\ldots,y^{(K)}$ and $z$.
Let us single out these giving the lowest order derivatives~:
\begin{equation}
\begin{array}{l}\label{uvij}
u^{(k-i_0)}=f^1(u,\ldots,u^{(k-1-i_0-1)},v,\ldots,v^{(\ell-j_0-1)},x,\ldots,x^{(K)},z,y,\ldots,y^{(K)}),\\
v^{(\ell-j_0)}=f^2(u,\ldots,u^{(k-1-i_0-1)},v,\ldots,v^{(\ell-j_0-1)},x,\ldots,x^{(K)},z,y,\ldots,y^{(K)}).
\end{array}
\end{equation}

Let us prove that, provided that $(x,y,z)$ is a solution of
(\ref{sys3}), system (\ref{uvij}) is equivalent
to (\ref{eqq101})-(\ref{eqq102})-(\ref{eqq103}), \ie\  to $\Gamma(u,v)=(x,y,z)$.
It is obvious that any $t\mapsto(u(t),v(t),x(t),y(t),z(t))$ that satisfies (\ref{sys3}),
(\ref{eqq101}), (\ref{eqq102}) and (\ref{eqq103}) also satisfies (\ref{uvij}),
because these equations were obtained from consequences of those. Conversely, let
$t\mapsto(u(t),v(t),x(t),y(t),z(t))$ be such that (\ref{sys3}) and
(\ref{uvij}) are satisfied; differentiating (\ref{uvij})
and substituting each time $\dot{z}$ from (\ref{sys3}) and
$(u^{(k-i_0)},v^{(\ell-j_0)})$ from (\ref{uvij}), one obtains
\begin{equation}\!
\begin{array}{l}\label{uvijDIFF}
u^{(k-i_0+i)}=f^{1,i}(u,\ldots,u^{(k-1-i_0-1)},v,\ldots,v^{(\ell-j_0-1)},x,\ldots,x^{(K+i)},z,y,\ldots,y^{(K+i)}),\
i\in\NN,\\
v^{(\ell-j_0+j)}=f^{2,j}(u,\ldots,u^{(k-1-i_0-1)},v,\ldots,v^{(\ell-j_0-1)},x,\ldots,x^{(K+j)},z,y,\ldots,y^{(K+j)}),\
j\in\NN.
\end{array}\!\!\!
\end{equation}
Now, substitute the values of $u^{(k-i_0)},\ldots,u^{(k)}$,
$v^{(\ell-j_0)},\ldots,v^{(\ell)}$ from (\ref{uvijDIFF}) into (\ref{eqq101}),
(\ref{eqq102}) and (\ref{eqq103}); either the obtained relations are
identically satisfied, and hence it is true that any solution of (\ref{sys3}) and
(\ref{uvij}) also satisfies (\ref{eqq101})-(\ref{eqq102})-(\ref{eqq103}), or
one obtains at least one relation of the form (recall that $k\leq\ell$):
$$
R(u,\ldots,u^{(k-1-i_0-1)},v,\ldots,v^{(\ell-j_0-1)},x,\ldots,x^{(K+\ell)},z,y,\ldots,y^{(K+\ell)})=0.
$$
This relation has been obtained (indirectly) by differentiating and combining 
(\ref{sys3})-(\ref{eqq101})-(\ref{eqq102})-(\ref{eqq103}).
This is absurd because
(\ref{eqq101})-(\ref{eqq102})-(\ref{eqq103})-(\ref{eqq102i})-(\ref{eqq102K}) are the only independent
relations of order $k,\ell$ obtained by differentiating and combining\footnote{
In other words,
(\ref{eqq101})-(\ref{eqq102})-(\ref{eqq103})-(\ref{eqq102i})-(\ref{eqq102K}),
as a system of ODEs in $u$ and $v$, is formally integrable (see e.g. \cite[Chapter IX]{Brya-Che-Gar-G-G91}). 
This means, for a systems of ODEs with independent variable $t$, that no new independent equation of the same orders ($k$ with
  respect to $u$ and $\ell$ with respect to $v$) can be obtained by
  differentiating 
  and combining these equations. 
  It is known \cite[Chapter IX]{Brya-Che-Gar-G-G91} that a sufficient
  condition is that this is true when differentiating only once and the system
  allows one to express the highest order derivatives as functions of the others.
  Formal integrability also means that, given any initial condition
  $(u(0),\ldots,u^{(k)}(0),v(0),\ldots,v^{(\ell)}(0))$ that satisfies these
  relations, there is a solution of the system of ODEs with these initial conditions.
}
(\ref{eqq101})-(\ref{eqq102})-(\ref{eqq103})
because, on the one hand, since $D^Kp\neq0$, differentiating more (\ref{eqq102K}) and (\ref{eqq101})
will produce higher order differential equations in which higher order
derivatives cannot be eliminated, and on the other hand, differentiating 
(\ref{eqq103}) and substituting $\dot{z}$ from (\ref{sys3}), $u^{(k)}$
from (\ref{eqq101}) and $\dot{y}$ from (\ref{eqq102i}) for $i=1$ yields the
trivial $0=0$ because $p$ is a solution of $\sysedp$, see the proof of
point~\ref{def2} above.

We have now established that, for $(x,y,z)$ a solution of (\ref{sys3}), $\Gamma(u,v)=(x,y,z)$ is equivalent to
(\ref{uvij}).
Using
Cauchy Lipschitz theorem with continuous dependence on the parameters, one can define a continuous map
$s:\mathcal{V}\to\mathcal{U}$ mapping a germ  $(x,y,z)$ to the unique
germ of solution of (\ref{uvij}) with fixed initial condition
$(u,\ldots,u^{(k-i_0-1)},v,\ldots,v^{(\ell-j_0-1)})=(u_0,\ldots,u_0^{(k-i_0-1)},v_0,\ldots,v_0^{(\ell-j_0-1)})$.
Then
$s$ is a continuous right inverse of $\Gamma$, \ie\  $\Gamma\circ s=Id$.
This proves point~\ref{def1}.
\end{proof}

\subsection{On (non-)existence of regular solutions of system $\sysedp$}
\label{sec-EDP}

\begin{cnjctr}
  \label{conj-EDP} For any real analytic functions $\gamma$ and $\delta$ (with
  $\gamma_4\neq0$), and any integers $k,\ell$, the partial differential system
  $\sysedp$ (see (\ref{EDPp})) does not admit any regular solution $p$.
\end{cnjctr}

An equivalent way of stating this conjecture is: ``the equations $ED^ip=0$,
for $1\leq i\leq k+\ell-2$, are consequences of (\ref{EDPp})''. Note that  ``$ED^ip=0$'' in
fact encodes several partial differential relations on $p$; see Remark~\ref{rmk-spur}.
If $\gamma$ and $\delta$ are polynomials, this can be easily
phrased in terms of the differential ideals in the set of polynomials with
respect to the variables $u,\ldots,u^{(k-1)},x,v,\ldots,v^{(\ell-1)}$ with
$k+\ell+1$ commuting derivatives (all the partial derivatives with respect to
these variables).

This is still a conjecture for general integers $k$ and $\ell$, but
we prove it for ``small enough'' $k,\ell$, namely
if one of them is smaller than 3 or if $k=\ell=3$.
The following statements assume $k\leq\ell$ (see remark~\ref{rmk-kl}).
\begin{prpstn}
\label{prop-edp}
If system $\sysedp$, with $k\leq\ell$, admits a regular solution, then $k\geq3$,
$\ell\geq4$ and the determinant 
\begin{equation}
  \label{eq:det2}\left|
    \begin{array}{lll}
p_{u^{(k-1)}}&p_{u^{(k-2)}}&p_{u^{(k-3)}} \\
p_{xu^{(k-1)}}&p_{xu^{(k-2)}}&p_{xu^{(k-3)}} \\
p_{xxu^{(k-1)}}&p_{xxu^{(k-2)}}&p_{xxu^{(k-3)}}
    \end{array}
\right|
\end{equation}
is a nonzero real analytic function.
\end{prpstn}
\begin{proof}
  Straightforward consequence of Lemma~\ref{lemAAA} and the three following
  lemmas, proved in appendix~\ref{app-sys}.
\end{proof}

\begin{lmm}
\label{lem-1}
If $p$ is a solution of system $\sysedp$
and either $k=1$ or 
$\left|
    \begin{array}{ll}
p_{u^{(k-1)}}&p_{u^{(k-2)}}\\
p_{xu^{(k-1)}}&p_{xu^{(k-2)}}
    \end{array}
\right|=0$, 
then around each point such that $p_{u^{(k-1)}}\neq0$, there exists a function
$\alpha$ of two variables such that a relation $p_x=\alpha(x,p)$ holds
identically on a neighborhood of that point.
\end{lmm}
\begin{lmm}
\label{lem-2}
Suppose that $p$ is a solution of $\sysedp$ with
\begin{equation}
  \label{inv2}
\ell\geq k\geq2\;,\ \ \ p_{u^{(k-1)}}\neq0\;,\ \ \ 
  \left|
    \begin{array}{ll}
p_{u^{(k-1)}}&p_{u^{(k-2)}}\\
p_{xu^{(k-1)}}&p_{xu^{(k-2)}}
    \end{array}
\right|
\neq0\ .
\end{equation}
If either $k=2$ or the determinant (\ref{eq:det2}) is identically zero, 
then, around any point where the two quantities in (\ref{inv2}) are nonzero,
there exists a function
$\alpha$ of three variables such that a relation $p_{x,x}=\alpha(x,p,p_x)$ holds
identically on a neighborhood of that point.
\end{lmm}
\begin{lmm}
\label{lem-3}
Let $k=\ell=3$.
For any solution $p$ of $\sysedpl\gamma\delta33$, in a neighborhood of any point where the
determinant (\ref{eq:det2}) is nonzero, there exist two functions $\alpha$
and $\psi$ of four variables such that
$p_{xxx}=\alpha(x,p,p_x,p_{xx})$ and
$Fp_{xx}+\tau p_{xxu^{(k-1)}}=\psi(x,p,p_x,p_{xx})$
identically on a neighborhood of that point.
\end{lmm}

\section{Remarks on the case where $S=T=0$.}
\label{sec-ST0}

\subsection{Geometric meaning of the differential form {\boldmath$\omega$ and the condition $S=T=0$}} 

For $(x,y,z)$ such that the set
$\Lambda=\{\lambda\in\RR,(x,y,z,\lambda)\in\Omega\}$ is nonempty,  
(\ref{sys3}) defines, by varying $\lambda$ in $\Lambda$ and
$\dot{x}$ in $\RR$, a surface $\Sigma$ in [the tangent space at
$(x,y,z)$ to] $\RR^3$. Fixing $\lambda$ in $\Lambda$ and varying
$\dot{x}$ in $\RR$ yields a straight line $S_\lambda$ (direction
$(1,z,g(x,y,z,\lambda))$). Obviously, $\Sigma=\bigcup_{\lambda\in\Lambda}S_\lambda$; 
$\Sigma$ is a ruled surface.
For each $\lambda\in\Lambda$, let $P_\lambda$ be the \emph{osculating
  hyperbolic paraboloid to $\Sigma$ along $S_\lambda$}, \ie\ the
unique\footnote{General hyperbolic paraboloid: 
  $\left(a^{11}\dot{x}+a^{12}Y+a^{13}Z\right) 
  \left(a^{21}\dot{x}+a^{22}Y+a^{23}Z\right)
  +a^{31}\dot{x}+a^{32}Y+a^{33}Z+a^0=0$, where the matrix
  $[a^{ij}]$ is invertible and $Y,Z$ stand for $\dot{y}-z\dot{x}-\lambda,
  \dot{z}-g\dot{x}-h$. It contains $S_\lambda$
  if and
  only if $a^{11}=a^{31}=a^0=0$. Contact at order 2 means $a^{13}=0$, $a^{33}=-a^{12}a^{21}/g_{4}$,
  $a^{32}=-h_{4}a^{33}$, $a^{22}=\frac12a^{21}(g_{4}h_{44}-g_{44}h_{4})/{g_{4}}^2$,
  $a^{23}=\frac12a^{21}g_{44}/{g_{4}}^2$.
  Normalization: $a^{12}=a^{21}=1$.
}
such quadric that contains $S_\lambda$ and has a contact of order 2 with
$\Sigma$ at \emph{all} points of $S_\lambda$. Its equation is
\begin{displaymath}
\textstyle\left( \dot{y}-z\dot{x}-\lambda \right)  
\left( \dot{x}+{\frac { h_{{44}}g_{{4}}-g_{{44}}h_{{4}} }{{2\,g_{{4}}}^{2}}}\left( \dot{y}-z\dot{x}-\lambda \right)+{\frac {g_{{44}}  }{2\,{g_{4}}^{2}}}\left( \dot{z}-g\dot{x}-h \right) \right) -{\frac {\dot{z}-g\dot{x}-h}{g_{{4}}}}+{\frac {h_{{4}}  }{g_{{4}}}}\left( \dot{y}-z\dot{x}-\lambda \right)=0
\end{displaymath}
where we omitted the argument $(x,y,z,\lambda)$ of $h$ and $g$. With $\omega$,
$\omega^1$, $\eta$ defined in (\ref{om})
and $\dot{\xi}$ the vector with coordinates $\dot{x},\dot{y},\dot{z}$, the
above equation reads
\begin{displaymath}
  -\,\left(\langle\omega^1,\dot{\xi}\rangle-\lambda\right)
  \frac{\langle\omega,\dot{\xi}\rangle + \left(
        h_{{44}}g_{{4}}-g_{{44}}h_{{4}} \right)\lambda + g_{{44}}h}
  {2{g_{4}}^{2}}
\,-\, \frac{ \langle\eta,\dot{\xi}\rangle-h}{g_{4}} =0\ ,
\end{displaymath}
that can in turn be rewritten 
$\langle\omega^1,\dot{\xi}\rangle\langle\omega,\dot{\xi}\rangle
-\langle\omega^3,\dot{\xi}\rangle-a^0=0$, with $\omega^3$ and $a^0$ some
differential form and function;
$\omega$, $\omega^3$ and $a^0$ are uniquely defined up to multiplication by a
non-vanishing function; they encode how the ``osculating hyperbolic paraboloid''
depends on $x,y,z$ and $\lambda$.

\smallskip
We will have to distinguish the case when $S$ and $T$, whose explicit
expressions derive from (\ref{om}) and (\ref{STJ}):
\begin{equation}\label{defST}
S=2\,g_{4}\,g_{4,4,4}\;-\;3\,{g_{4,4}}^2\;,\ \ \ 
T=2\,g_{4}\,h_{4,4,4}\;-\;3\,g_{4,4}\,h_{4,4}\;,
\end{equation}
are zero.
From (\ref{STJ}), it means that the Lie derivative of
  $\omega$ along $\partial/\partial\lambda$ is co-linear to $\omega$, and this
is classically equivalent to a decomposition $\omega=k\,\hat{\omega}^2$ where
$k\neq0$ is a function of the four variables $x,y,x,\lambda$ but 
$\hat{\omega}^2$ is a differential form in the \emph{three} variables $x,y,z$, the
first integrals of $\partial/\partial\lambda$. Then, 
one can prove that the form $\hat{\omega}^3=\omega^3/k$ and the function $\hat{a}^0=a^0/k$
also involve the variables $x,y,z$ only.
From $\omega$'s expression, one can
take for instance $k=g_{44}$ or $k=gg_{44}-2{g_{4}}^2$ (they do not vanish
simultaneously because $g_{4}$ does not vanish).
Hence $S=T=0$ if and only if, for each fixed $(x,y,z)$, the osculating
hyperbolic paraboloid $P_\lambda$ in fact does not depend on $\lambda$ \ie\
the surface $\Sigma$ itself is a hyperbolic paraboloid, its
equation being
\begin{equation}
  \label{eqq50}
  \langle\omega^1,\dot{\xi}\rangle
\langle\hat{\omega}^2,\dot{\xi}\rangle
+
\langle\hat{\omega}^3,\dot{\xi}\rangle
+
\hat{a}^0=0\ ,
\end{equation}
where $\dot\xi$ is the vector of coordinates $\dot{x},\dot{y},\dot{z}$.
This yields the following proposition\footnote{ 
  We introduced the osculating hyperbolic paraboloid because it gives some
  geometric insight on $\omega$, $S$ and $T$, but it is not formally
  \emph{needed}: Proposition~\ref{propST0} can be
  stated without it, and proved as
  follows, based on (\ref{defST}) (see also \cite{Avan05th}): the general
  solution of $S=0$ is a 
  linear fractional expression $\displaystyle g= 
    \bigl.\bigl(\hat{b}^0+\hat{b}^1\lambda \bigr)\bigr/
    \bigl(\hat{c}^0+\hat{c}^1\lambda\bigr)$ where $\hat{b}^0$, $\hat{b}^1$,
    $\hat{c}^0$, $\hat{c}^1$ are functions of $x,y,z$ only
  ---this is known, for $S/(g_{4})^2$ is the Schwartzian derivative of $g$
  with respect to its 4\textsuperscript{th} argument, but anyway elementary---
   and $g_4\neq0$ translates into $\hat{b}^0\hat{c}^1-\hat{b}^1\hat{c}^0\neq
   0$; then $T=0$ yields $\displaystyle h= 
    \bigl.\bigl(\hat{a}^0+\hat{a}^1\lambda+\hat{a}^2\lambda^2 \bigr)\bigr/
    \bigl(\hat{c}^0+\hat{c}^1\lambda\bigr)$
  with $\hat{a}^0$, $\hat{a}^1$, $\hat{a}^2$ functions of $x,y,z$.
  With such $g$ and $h$, multiplying both sides of (\ref{sys3}) by
  $\hat{c}^0+\hat{c}^1\lambda$ yields the equation in
  Proposition~\ref{propST0}. 
}, 
where the functions $\hat{a}^0$, $\hat{a}^1$,
$\hat{a}^2$, $\hat{b}^0$, $\hat{b}^1$, $\hat{c}^0$, $\hat{c}^1$ of $x,y,z$ are
defined by
\begin{equation}
  \label{om12}
\hat{\omega}^2=\frac{\omega}{k}=\hat{b}^1\xdif{x}+\hat{a}^2\omega^1-\hat{c}^1\xdif{z}
,\ \
\hat{\omega}^3=\frac{\omega^3}{k}=\hat{b}^0\xdif{x}+\hat{a}^1\omega^1-\hat{c}^0\xdif{z}
,\ \
\hat{a}^0=\frac{a^0}{k}\;.
\end{equation}

\begin{prpstn}
\label{propST0}
If $S$ and
$T$, given by (\ref{om})-(\ref{STJ}), or (\ref{defST}), are identically zero
on $\Omega$, 
then, for any $(x_0,y_0,z_0,\lambda_0)$ in $\Omega$, there exist an open set
$W\subset\RR^3$, an open interval $I\subset\RR$, with
$(x_0,y_0,z_0,\lambda_0)\in W\times I\subset\Omega$, and seven smooth functions $W\to\RR$ denoted by $\hat{a}^0$, $\hat{a}^1$, $\hat{a}^2$,
$\hat{b}^0$, $\hat{b}^1$, $\hat{c}^0$, $\hat{c}^1$ such that
$\hat{c}^0+\hat{c}^1\lambda$ does not vanish on $W\times I$,
$\hat{c}^1\hat{b}^0-\hat{b}^1\hat{c}^0$ does not vanish on $W$,
 and, for $(x,y,z,\lambda)\in W\times I$, $\dot{x}\in\RR$ and $\dot{z}\in\RR$,
 equation (\ref{sys3}) is equivalent to 
\begin{displaymath}
\lambda\left(\hat{b}^1\!(\!x,\!y,\!z\!)\dot{x}+\hat{a}^2\!(\!x,\!y,\!z\!)\lambda
  -\hat{c}^1\!(\!x,\!y,\!z\!)\dot{z} 
\right)
+
\left( \hat{b}^0\!(\!x,\!y,\!z\!)\dot{x}+\hat{a}^1\!(\!x,\!y,\!z\!)\lambda -\hat{c}^0\!(\!x,\!y,\!z\!)\dot{z}\right)
+
\hat{a}^0\!(\!x,\!y,\!z\!)=0.
\end{displaymath}
\end{prpstn}

\subsection{A parameterization of order $(1,2)$ if $S=T=J=0$}

It is known~\cite{Pome97cocv} that system (\ref{sys3})
is $(x,u)$-flat (see section~\ref{sec-flat}) if $S=T=J=0$. For the sake of completeness, let re-state this result
in terms of parameterization. We start with the following  particular case of (\ref{sys3}):
\begin{equation}
  \label{eq:new2}
  \dot{z} = \kappa(x,y,z) \,\dot{x} \,\lambda
+a(x,y,z) \,\lambda
+b(x,y,z) \,\dot{x}+c(x,y,z)
\ \ \mbox{with}\ \ 
\lambda=\dot{y}-z\dot{x}
\end{equation}
where $\kappa$ does not vanish on the domain where it is defined.
Note that Example~\ref{ex-xu} was of this type with $\kappa=1$, $a=b=0$, $c=y$.
For short, define the following vector fields:
$$
X^0=c\partialx{z}\,,\ 
X^1=\partialx{x}+z\partialx{y}+b\partialx{z}\,,\ 
X^2=\partialx{y}+a\partialx{z}\,,\ 
X^3=\kappa\partialx{z}\ .
$$
Note that, for $h$  an arbitrary smooth function of
$x$, $y$ and $z$, $X^0h$, $X^1h$, $X^2h$, $X^3h$ also depend on $x,y,z$ only.

\begin{lmm} 
\label{lmm-canonxulin}
System (\ref{eq:new2}) admits a parameterization of order (1,2) at
  any $(x_0,y_0,z_0,\dot{x}_0,\dot{y}_0$, $\ddot{x}_0,\ddot{y}_0)$ such that
  \begin{equation}
    \label{eq:new4}
    \kappa\,\ddot{x}_0+\kappa^2{\dot{x}_0}^3
    +\bigl(X^1\kappa-X^3b+2a\kappa\bigr){\dot{x}_0}^2
    +\bigl(X^1a+X^0\kappa-X^3c-X^2b+a^2\bigr){\dot{x}_0} +X^0a-X^3c\neq0.
  \end{equation}
\end{lmm}
\begin{proof}
  From (\ref{eq:new4}), the two vector fields $Y=X^2+\dot{x}X^3$ and 
$Z=[\,X^0+\dot{x}X^1\,,\,X^2+\dot{x}X^3\,]+\ddot{x}X^3$ are linearly
independent at point $(x_0,y_0,z_0,\dot{x}_0,\ddot{x}_0)$.
Let then $h$ be a function of $(x,y,z,\dot{x})$ such that $Yh=0$ 
and $Zh\neq 0$; its ``time-derivative along system
(\ref{eq:new2})'', given by
$\dot{h}= X^0h +\bigl(X^1h\bigr)\dot{x} +\bigl(Yh\bigr)\lambda+\bigl(\partial h/\partial\dot{x})\ddot{x}$,
does not depend on $\lambda$: it is a function of
$(x,y,z,\dot{x},\ddot{x})$; also, since $Yh=0$, one has $Y\dot{h}=Zh$;
finally, $Zh\neq 0$ implies that
$\xdif{h}\wedge\xdif{\dot{h}}\wedge\xdif{x}\wedge\xdif{\dot{x}}\wedge\xdif{\ddot{x}}\neq0$.
In turn, this implies that
$(x,y,z,\dot{x},\ddot{x})\mapsto(h(x,y,z,\dot{x}),\dot{h}(x,y,z,\dot{x},\ddot{x}),x,\dot{x},\ddot{x})$
defines a local diffeomorphism at $(x_0,y_0,z_0,\dot{x}_0,\ddot{x}_0)$. Let $\psi$ and $\chi$ be the two functions of
five variables such that the inverse of that local diffeomorphism is
$(u,\dot{u},v,\dot{v},\ddot{v})\mapsto
(v,\psi(u,\dot{u},v,\dot{v},\ddot{v}),\chi(u,\dot{u},v,\dot{v},\ddot{v}),\dot{v},\ddot{v})$.
The parameterization (\ref{para1}) is given by~:
$x=v,y=\psi(u,\dot{u},v,\dot{v},\ddot{v}),z=\chi(u,\dot{u},v,\dot{v},\ddot{v})$.
\end{proof}

\begin{thrm}
  \label{prop-plats}
If $S=T=J=0$, then system (\ref{sys3}) admits a parameterization of order $(1,2)$ at any 
$(x_0,y_0,z_0,\dot{x}_0$, $\dot{y}_0,\ddot{x}_0,\ddot{y}_0) \in
\bigl(\widehat{\Omega}\times\RR^2\bigr)\setminus F$,
where $F\subset\widehat{\Omega}\times\RR^2$ is closed with empty interior. 
\end{thrm}
\begin{proof}
  From Proposition~\ref{propST0}, (\ref{sys3}) and (\ref{eqq50}) are
  identical. Since 
  $\xdif{\hat{\omega}^2}\wedge\hat{\omega}^2=0$ (see (\ref{STJ})-(\ref{om12})),
  there is a local
  change of coordinates $(\tilde{x},\tilde{y},\tilde{z})=P(x,y,z)$ such that
  $\hat{\omega}^2=k'\xdif{\tilde{x}}$ and
  $\omega^1=k''\big(\xdif{\tilde{y}}-\tilde{z}\xdif{\tilde{x}}\big)$ with
  $k'\neq0$, $k''\neq0$. Hence $P$ transforms (\ref{eqq50}) into
  (\ref{eq:new2}), for some $\kappa,a,b,c$.
  Lemma~\ref{lmm-canonxulin} gives $\varphi,\psi,\chi$ defining a
  parameterization of order (1,2) for this system. Then
  $P^{-1}\circ\varphi,P^{-1}\circ\psi,P^{-1}\circ\chi$ define one for
  the original system (\ref{sys3}), or (\ref{eqq50});
  $\bigl(\widehat{\Omega}\times\RR^2\bigr)\setminus F$ is the inverse image by
  $P$ of the set defined by (\ref{eq:new4}).
\end{proof}

\subsection{A normal form if $S=T=0$ and $J\neq0$}

\begin{prpstn}
\label{lem11}
Assume that the functions $g$ and $h$ defining system (\ref{sys3}) are such
that $S$ and $T$ defined by (\ref{STJ}) or (\ref{defST}) are identically zero
on $\Omega$, and let
$(x_0,y_0,z_0,\lambda_0)\in\Omega$ be such that
$J(x_0,y_0,z_0,\lambda_0)\neq0$.

There exist  an open set
$W\subset\RR^3$ and an open interval $I\subset\RR$ such that
$(x_0,y_0,z_0,\lambda_0)\in W\times I\subset\Omega$, a smooth
diffeomorphism $P$ from $W$ to $P(W)\subset\RR^3$ and six smooth functions
$P(W)\to\RR$ denoted $\kappa$, $\alpha$,
$\beta$, $a$, $b$, $c$ such that, with the 
change of coordinates
$(\tilde{x},\tilde{y},\tilde{z})=P(x,y,z)$,  system (\ref{sys3}) 
reads
\begin{equation}
  \label{eq:bof}
  \dot{\tilde{z}}=\kappa(\tilde{x},\tilde{y},\tilde{z})\left(\dot{\tilde{y}}-\alpha(\tilde{x},\tilde{y},\tilde{z})\,\dot{\tilde{x}}\right)
\left(\dot{\tilde{y}}-\beta(\tilde{x},\tilde{y},\tilde{z})\,\dot{\tilde{x}}\right)
\;+\;a(\tilde{x},\tilde{y},\tilde{z})\,\dot{\tilde{x}}+b(\tilde{x},\tilde{y},\tilde{z})\,\dot{\tilde{y}}+c(\tilde{x},\tilde{y},\tilde{z})
\end{equation}
and none of the functions
 $\kappa$, $\alpha-\beta$,
$\alpha_3$ and $\beta_3$ vanish on $W$.
\end{prpstn}

\begin{proof}
From Lemma~\ref{propST0}, we consider system (\ref{eqq50}).  
Let $P^1,P^2$ be a pair of independent first integrals of the vector field
$\hat{c}^1\left(\partialx{x}+z\partialx{y}\right) +\hat{b}^1\partialx{z}$;
from (\ref{om12}), $\omega^1,\hat{\omega}^2$ span the annihilator of this
vector field, and hence are independent linear combinations of $\xdif{P}^1$ and
$\xdif{P}^2$:
possibly interchanging $P^1$ and $P^2$ or adding one to the other,
there exist smooth functions
$k^1,k^2,f^1,f^2$ such that
$\hat{\omega}^i=k^i\left(\xdif{P}^2-f^i\xdif{P}^1\right)$, $f^1-f^2\neq0$, $k^i\neq0$,
$i=1,2$.
Now, take for $P^3$ any function such that
$\xdif{P}^1\wedge\xdif{P}^2\wedge\xdif{P}^3\neq0$;
decomposing $\hat{\omega}^3$, we get three smooth
functions $ p^0, p^1, p^2$ such that $\hat{\omega}^3\ =\ p^0\left(-\xdif{P}^3+
p^1\xdif{P}^1+ p^2\xdif{P}^2\right)$, $p^0\neq0$.
The change of coordinates
$P=(P^1,P^2,P^3)$ does transform system
(\ref{eqq50}) into (\ref{eq:bof}) with
$$
\kappa=\frac{k^1k^2}{ p^0}\circ P^{-1},\ \ 
\alpha=f^1\circ P^{-1},\ \ \beta=f^2\circ P^{-1},\ \ a= p^1\circ P^{-1},\ \
b= p^2\circ P^{-1},\ \ c=\frac{\hat{a}^0}{ p^0}\circ P^{-1}\;.
$$
$\kappa$ and $\alpha-\beta$ are nonzero because $f^1-f^2$, $k^1$ and $k^2$
are.
$\alpha_3$ and $\beta_3$ are nonzero because the inverse images of 
$\alpha_3\xdif{\tilde{x}}\wedge\xdif{\tilde{y}}\wedge\xdif{\tilde{z}}$ and
$\beta_3\xdif{\tilde{x}}\wedge\xdif{\tilde{y}}\wedge\xdif{\tilde{z}}$ by $P$
are $\xdif{P}^1\wedge\xdif{P}^2\wedge\xdif{f}^i$ for $i=1,2$, that are equal,
by construction, to $\xdif{\omega}^1\wedge\omega^1/(k^1)^2$ and
$\xdif{\hat{\omega}}^{2}\wedge\hat{\omega}^2/(k^2)^2$, which are both nonzero
(the second one because $J\neq0$).
\end{proof}

Note that (\ref{eq:bof}) is not in the form (\ref{sys3})
unless $\alpha=\tilde{z}$ or $\beta=\tilde{z}$. This suggests, since $\alpha_3\neq0$ and $\beta_3\neq0$, the following
local changes of coordinates $A$ and $B$,  that both turn
(\ref{eq:bof}) to a new system of the form (\ref{sys3}):
\begin{equation}
  \label{diffST0}
  (\tilde{x},\tilde{y},\tilde{z})\mapsto A(\tilde{x},\tilde{y},\tilde{z})=(\tilde{x},\tilde{y},\alpha(\tilde{x},\tilde{y},\tilde{z}))\ \ \ \mbox{and}\ \ \
  (\tilde{x},\tilde{y},\tilde{z})\mapsto B(\tilde{x},\tilde{y},\tilde{z})=(\tilde{x},\tilde{y},\beta(\tilde{x},\tilde{y},\tilde{z}))\;.
\end{equation}
These two systems of the form (\ref{sys3}) correspond to two choices $h^1,g^1$
and $h^2,g^2$ instead of the original $h,g$, and they yield, according to (\ref{defgamma}) and
(\ref{defdelta}), two possible sets of
functions $\gamma$ and $\delta$. These will be used in Theorem~\ref{th0};
let us give their explicit expression~: 
\begin{equation}
  \label{case1}
  \gamma^i(x,y,z,w)=
\frac{w-m^{i,0}(x,y,z)}{m^{i,1}(x,y,z)}
\ ,\ \ 
\delta^i=n^{i,0}+n^{i,1}\gamma+n^{i,2}\gamma^2
\ ,\ \ i\in\{1,2\}
\end{equation}
with (these are obtained from each other by interchanging $\alpha$ and $\beta$)~:
\begin{equation}
  \label{case2}
  \begin{array}{l}
m^{1,0}=\left(\alpha_1+\alpha\,\alpha_2+(a+b\,\alpha)\,\alpha_3\right)\circ A^{-1}\,,\ \ 
m^{1,1}=\left(\kappa\,\alpha_3\,(\alpha-\beta)\right)\circ A^{-1}\,,\\
n^{1,0}=\alpha_3\circ A^{-1}\,,\ \ 
n^{1,1}=(\alpha_2+b\,\alpha_3)\circ A^{-1}\,,\ \ 
n^{1,2}=(\kappa\,\alpha_3)\circ A^{-1}\,,
\\[1ex]
m^{2,0}=\left(\beta_1+\beta\,\beta_2+(a+b\,\beta)\,\beta_3\right)\circ B^{-1}\,,\ \ 
m^{2,1}=\left(\kappa\,\beta_3\,(\beta-\alpha)\right)\circ B^{-1}\,,\\
n^{2,0}=\beta_3\circ B^{-1}\,,\ \ 
n^{2,1}=(\beta_2+b\,\beta_3)\circ B^{-1}\,,\ \ 
n^{2,2}=(\kappa\,\beta_3)\circ B^{-1}\,.
  \end{array}
\end{equation}
\begin{xmpl}
  \label{ex:2gamma}
System (\ref{eq:ex}-b) in Example~\ref{3ex} is already as in (\ref{eq:bof}). The above choices are, for this system:
\begin{equation}
  \label{eq:ex2gamma}
    \gamma^1(x,y,z,w)=w
,\ 
\delta^1(x,y,z,w)=y+w^2
,\ \gamma^2(x,y,z,w)=-w
,\ 
\delta^2(x,y,z,w)=y+w^2.
\end{equation}
\end{xmpl}

\section{Main results}
\label{sec-main}

We gather here our main results in a synthetic manner.
They rely on precise \emph{local} results from other sections~: 
sufficient (sections~\ref{sec-ST0} and \ref{sec-edpsuff})
or necessary (section~\ref{sec-nec}) conditions for parameterizability,
results on solutions of the partial differential system $\sysedp$
(section~\ref{sec-EDP}) and on the relation between flatness and
parameterizability (section~\ref{sec-flat}).
We are not able to give local precise necessary and sufficient conditions at a
given point (jet) because singularities are not the same for necessary and for
sufficient conditions; instead, we use the ``somewhere'' notion as in
Definitions~\ref{def-Kreg-gen} and \ref{def-param-bof}.

\begin{thrm}\label{th-cns}
System (\ref{sys3}) admits a parameterization of order $(k,\ell)$ somewhere in $\Omega$ if and only if 
\begin{enumerate}
\item  \label{cns1}either $S=T=J=0$ on $\Omega$ (in this case, one can take
  $(k,\ell)=(1,2)$),
\item \label{cns2}or $S=T=0$ on $\Omega$ and one of the two systems
  $\sysedpl{\gamma^1}{\delta^1}k\ell$ or $\sysedpl{\gamma^2}{\delta^2}k\ell$ 
  with $\gamma^i$, $\delta^i$ given by (\ref{case1})-(\ref{case2}), admits
  a regular solution somewhere in $\widehat{\Omega}$.
\item \label{cns3}or $S$ and $T$ are not both identically zero, and the system $\sysedp$
  with $\gamma$ and $\delta$ defined from $g$ and $h$ according to
  (\ref{defgamma}) and (\ref{defdelta}) admits
  a regular solution somewhere in $\widehat{\Omega}$.
\end{enumerate}
\end{thrm}
\begin{proof}
Sufficiency~:
the parameterization is provided, away from an explicitly described set of
singularities, by Theorem~\ref{prop-plats} if
point~\ref{cns1} holds, and by Theorem~\ref{edpparam} if one of the two
other points holds.
For necessity, assume that there is a
parameterization of order $(k,\ell)$ at a point $(x,y,z,\dot{x},\dot{y},\ldots,x^{(L)},y^{(L)})$
in $\bigl(\widehat{\Omega}\times\RR^{2L-2}\bigr)\backslash F$. From Theorems
\ref{thSTnot0} and \ref{th0}, it implies that one of the three points holds.
\end{proof}

\begin{xmpl}\label{ex:mainth}
  Consider again systems (a), (b) and (c) in (\ref{eq:ex}).
  From point~\ref{cns1} of the theorem, system (a) admits a parameterization
  of order (1,2), see also Example~\ref{ex-xu}. System (b) is concerned by
  point~\ref{cns2} of the theorem: it has a parameterization of order $k,\ell$
  if and only one of the two systems of PDEs 
\begin{equation}
\label{EDPp2}
\begin{array}{ll}
p_{u^{(k-1)}}\bigl(Fp_x-p-{p_{xx}}^2)\bigr)
-p_{xu^{(k-1)}}\bigl(Fp\pm p_{xx})\bigr)=p_{u^{(k-1)}}\,p_{xv^{(\ell-1)}}-p_{xu^{(k-1)}}\,p_{v^{(\ell-1)}}=0\;,
\\
p_{u^{(k-1)}}\neq 0\;,\ \ \ 
p_{v^{(\ell-1)}}\neq 0\;,\ \ \ 
p+{p_{xx}}^2\pm p_{xxx}\neq0
\end{array}
\end{equation}
admits a ``regular solution''. Point~\ref{cns3} of the theorem is relevant to
system (c) because $S\neq0$: (c) admits a parameterization of
order $k,\ell$ if and only there is a ``regular solution'' $p$ to
\begin{equation}
\label{EDPp3}
\begin{array}{ll}
p_{u^{(k-1)}}\bigl(Fp_x-p)\bigr)
-p_{xu^{(k-1)}}\bigl(Fp-\sqrt{p_{xx}})\bigr)=p_{u^{(k-1)}}\,p_{xv^{(\ell-1)}}-p_{xu^{(k-1)}}\,p_{v^{(\ell-1)}}=0\;,
\\
p_{u^{(k-1)}}\neq 0\;,\ \ \ 
p_{v^{(\ell-1)}}\neq 0\;,\ \ \ 
p-{p_{xxx}}/{2\sqrt{p_{xx}}} \neq0\ .
\end{array}
\end{equation}
If Conjecture~\ref{conj-EDP} is true, neither system (b) nor system (c) admits
a parameterization of any order.
\end{xmpl}

Theorem~\ref{th-cns} gives a central role to the system of PDEs $\sysedp$. It
makes Conjecture~\ref{conj-EDP} equivalent to Conjecture~\ref{conj-param}
below. 
Theorem~\ref{th-3} states that the conjecture is true for $k,\ell$ ``small enough''.
\begin{cnjctr}
\label{conj-param}
If $\xdif{\omega}\wedge\omega$ (or $(S,T,J)$) is not identically zero on
$\Omega$, then system (\ref{sys3}) does not admit a parameterization of any
order at any point (jet of any order).
\end{cnjctr}
\begin{thrm}\label{th-3}
If system (\ref{sys3}) admits a parameterization of order $(k,\ell)$, with $k\leq\ell$, at some
jet, then either $S=T=J=0$ or $k\geq3$ and $\ell\geq4$.
\end{thrm}
\begin{proof}
  This is a simple consequence of Theorem~\ref{th-cns} and
  Proposition~\ref{prop-edp}.
\end{proof}

\begin{rmrk}
\label{rmk-suff}
  If our Conjecture~\ref{conj-EDP} is correct, the systems $\sysedp$ never
  have any regular solutions, and the sufficiency part of Theorem~\ref{th-cns}
  (apart from case~\ref{cns1}) is essentially void, and so is
  Theorem~\ref{edpparam}. However, Conjecture~\ref{conj-EDP} is still a
  conjecture, and the interest of the sufficient conditions above is to
  make this conjecture, that only deals with a set of
  partial differential equalities and inequalities, equivalent to
  Conjecture~\ref{conj-param} below.
  For instance, if one comes up with a regular solution of some of these
  systems  $\sysedp$, this will yield a new class of systems that admit a
  parameterization.
\end{rmrk}

\begin{rmrk}[on recovering the results of \cite{Pome97cocv}]
\label{rmrk-old}
The main result in that reference can be phrased~: 
\begin{quote}
  ``~(\ref{sys4}) is $(x,u)$-dynamic linearizable (\ie\  $(x,u)$-flat) if
  and only if $S=T=J=0$~''~.
\end{quote}
Sufficiency is elementary in~\cite{Pome97cocv}; Theorem~\ref{prop-plats} implies it.
The difficult part is to prove that $S=T=J=0$ is necessary; that
proof is very technical in~\cite{Pome97cocv}: it relies on some simplifications performed via
computer algebra.
From our Proposition~\ref{prop-xu33}, $(x,u)$-flatness implies existence of a parameterization of
some order $(k,\ell)$ with $k\leq 3$ and $\ell\leq3$. Hence Theorem~\ref{th-cns}
does imply the above statement.
\end{rmrk}

\section{Necessary conditions}
\label{sec-nec}

\subsection{The case where $S$ and $T$ are not both zero}

The following lemma is needed to state the theorem.

\begin{lmm}\label{lemphi}
If $(S,T,J)\neq(0,0,0)$ and system (\ref{sys3}) admits a
parameterization $(\varphi,\psi,\chi)$ of order $(k,\ell)$ 
at point
$(x_0,y_0,z_0,\ldots,x^{(L)}_0,y^{(L)}_0)\in\RR^{2L+3}$,
then $\varphi_{u^{(k)}}$ is a nonzero real analytic function.
\end{lmm}
\begin{proof}
Assume a parameterization where $\varphi$ does not depend on $u^{(k)}$. Substituting in (\ref{sys3}) yields
$$
\dot{\chi}=h(\varphi,\psi,\chi,\dot{\psi}-\chi\dot{\varphi})
+g(\varphi,\psi,\chi,\dot{\psi}-\chi\dot{\varphi})\dot{\varphi}\ .
$$
Since $\dot{\varphi}$ does
not depend on $u^{(k+1)}$, differentiating twice with respect to $u^{(k+1)}$ yields
$$
\chi_{u^{(k)}}=\psi_{u^{(k)}}(h_{4}+g_{4}\dot{\varphi})\ ,\ \ \ 
0={\psi_{u^{(k)}}}^2(h_{4,4}+g_{4,4}\dot{\varphi})\ .
$$
If $\psi_{u^{(k)}}$ was zero, then, from the first relation,
$\chi_{u^{(k)}}$ would too, and this would contradict point~\ref{def3} in
Definition~\ref{defparam}; hence the second relation implies that
$h_{4,4}+g_{4,4}\dot{\varphi}$ is identically zero. From point~\ref{def1} in the
same definition, it implies that \emph{all} solutions of (\ref{sys3}) satisfy
the relation~:
$h_{4,4}(x,y,z,\dot{y}-z\dot{x})+g_{4,4}(x,y,z,\dot{y}-z\dot{x})\dot{x}=0$.
From Lemma~\ref{lem-absurd}, this implies that $h_{4,4}$ and
$g_{4,4}$ are the zero function of four variables, and hence $S=T=J=0$. This
proves the lemma.
\end{proof}

\begin{thrm}\label{thSTnot0}
Assume that either $S$ or $T$ is not identically zero on $\Omega$, and
that system (\ref{sys3}) admits a
parameterization of order $(k,\ell)$ at
$\mathcal{X}=(x_0,y_0,z_0,\dot{x}_0,\dot{y}_0,\ldots,x^{(L)}_0,y^{(L)}_0)\in\widehat{\Omega}\times\RR^{2L-2}$,
with $k,\ell,L$ some integers and 
$\varphi,\psi,\chi$ defined on $U\subset\RR^{k+\ell+2}$.

Then $k\geq1$, $\ell\geq 1$ and,  for any point
$(u_0,\ldots,u_0^{(k)},v_0,\ldots,v_0^{(\ell)})\in U$ (not necessarily sent to
$\mathcal{X}$ by the parameterization) such that
$$\varphi_{u^{(k)}}(u_0,\ldots,u_0^{(k)},v_0,\ldots,v_0^{(\ell)})\neq 0,$$ 
there exist a neighborhood $\ouvedp$ of 
$(u_0,\ldots,u^{(k-1)}_0,\,\varphi(u_0\,\cdots\,v_0^{(\ell)})\,,v_0,\ldots,v^{(\ell-1)}_0)$
in $\RR^{k+\ell+1}$ and a regular solution $p:\ouvedp\to\RR$ of $\sysedp$,
related to $\varphi,\psi,\chi$ by
(\ref{defphi}), (\ref{defpsi}) and (\ref{defchi}), the functions $\gamma$ and
$\delta$ being related to $g$ and
$h$ by (\ref{defgamma}) and (\ref{defdelta}).
\end{thrm}
\begin{rmrk}
\label{rmk-Lmin}
  The regular solution $p$ is $K$-regular for some positive integer $K\leq
  k+\ell-2$. If $L>K$, Theorem~\ref{edpparam} implies,
  possibly away from some singular values of
  $(x_0,y_0,z_0,\dot{x}_0,\dot{y}_0,\ldots$, $x^{(K)}_0,y^{(K)}_0)$, that system
  (\ref{sys3}) also admits a parameterization of order $(k,\ell)$ at
  $(x_0,y_0,z_0,\dot{x}_0,\dot{y}_0$, $\ldots,x^{(K)}_0,y^{(K)}_0)$. 
  See also Remark~\ref{rmrk-L}.
\end{rmrk}

\begin{proof}
Assume that system (\ref{sys3}) admits a parameterization
$(\varphi,\psi,\chi)$ of order $(k,\ell)$ at
$(x_0,y_0$, $z_0, \dot{x}_0,$ $\dot{y}_0$, $\ldots, x^{(L)}_0$, $y^{(L)}_0)$.
Since $\varphi_{u^{(k)}}$ does not vanish, one can apply the inverse function theorem to
the map
$$(u,\dot{u},\ldots,u^{(k)},v,\dot{v},\ldots,v^{({\ell})}) \mapsto
(u,\ldots,u^{(k-1)},\varphi(u,\ldots,u^{(k)},v,\ldots,v^{(\ell)}),v,\ldots,v^{({\ell})})$$
 and define locally a function $r$ of $k+\ell+2$ variables such that
 \begin{equation}
   \label{defr}
   \varphi(u,\dot{u},\ldots,u^{(k)},v,\dot{v},\ldots,v^{({\ell})})=x
\ \Leftrightarrow\ 
r(u,\dot{u},\ldots,u^{(k-1)},x,v,\dot{v},\ldots,v^{({\ell})})=u^{(k)}
\ .
 \end{equation}
Defining two functions $p,q$ by substitution of $u^{(k)}$ in $\psi$, $\chi$,
 the parameterization can be re-written implicitly as
\begin{equation} \label{eqp} 
\left\{\begin{array}{ll}
y=p(u,\dot{u},\ldots,u^{(k-1)},x,v,\dot{v},\ldots,v^{({\ell})}),\\
z=q(u,\dot{u},\ldots,u^{(k-1)},x,v,\dot{v},\ldots,v^{({\ell})}),\\
u^{(k)}=r(u,\dot{u},\ldots,u^{(k-1)},x,v,\dot{v},\ldots,v^{({\ell})}).
\end{array}\right. \end{equation}
We now work with this form of the parameterization and
$u,\dot{u},\ldots,u^{(k-1)},x,\dot{x},\ddot{x},\ldots\;v,\dot{v},\ldots$, $v^{({\ell})},v^{({\ell+1})},\ldots\;$
instead of 
$u,\dot{u},\ldots,u^{(k-1)},u^{(k)},\ldots\;v,\dot{v},\ldots,v^{({\ell})},v^{({\ell+1})},\ldots\;$.
In order to simplify notations, 
let us agree that, if $k=0$, the list $u,\dot{u},\ldots,u^{(k-1)}$ is
empty and any term involving the index $k-1$ is zero (same with $\ell-1$ if $\ell=0$).
Let us also define $\mathcal{P}$ and $\mathcal{Q}$ by
\begin{equation}
\label{PQ}
  \mathcal{P}=Fp+rp_{u^{(k-1)}}+v^{(\ell)}p_{v^{(\ell-1)}}+v^{(\ell+1)}p_{v^{(\ell)}}\,,\ \ 
  \mathcal{Q}=Fq+rq_{u^{(k-1)}}+v^{(\ell)}q_{v^{(\ell-1)}}+v^{(\ell+1)}q_{v^{(\ell)}}\ ,
\end{equation}
with $F$ given by (\ref{F}).
$\mathcal{P}$ and $\mathcal{Q}$ depend on
$u,\dot{u},\ldots,u^{(k-1)},x,v,\dot{v},\ldots,v^{({\ell+1})}$
but not on $\dot{x}$; $Fp$ and $Fq$ depend neither on $\dot{x}$ nor on $v^{({\ell+1})}$.
When substituting (\ref{eqp}) in (\ref{sys3}), using
$\dot{y}={\mathcal{P}}+\dot{x}p_x$ and $\dot{z}={\mathcal{Q}}+\dot{x}q_x$, 
one obtains~:
\begin{equation}\label{eq2}
{\mathcal{Q}}+\dot{x}q_x=h(x,p,q,\lambda)+g(x,p,q,\lambda)\dot{x}\ \ \
\mbox{with}\ \  \lambda={\mathcal{P}}+\dot{x}(p_x-q).
\end{equation}

Differentiating each side three times with respect to
$\dot{x}$, one obtains~:
\begin{eqnarray}
\label{eq201}
&&\hspace{-3.6em}
q_x=
\left(h_{4}(x,p,q,\lambda)+g_{4}(x,p,q,\lambda)\dot{x}\right)(p_x-q)+g(x,p,q,\lambda),\\
\label{eq202}
&&\hspace{-3.6em}
0=
\left(h_{4,4}(x,p,q,\lambda)+g_{4,4}(x,p,q,\lambda)\dot{x}\right)(p_x-q)^2+2g_{4}(x,p,q,\lambda)(p_x-q),\\
\label{eq203}
&&\hspace{-3.6em}
0=
\left(h_{4,4,4}(x,p,q,{\lambda})+g_{4,4,4}(x,p,q,{\lambda})\dot{x}\right)(p_x-q)^3+3g_{4,4}(x,p,q,\lambda)(p_x-q)^2.
\end{eqnarray}
Combining (\ref{eq202}) and (\ref{eq203}) to cancel the first term in each
equation, one obtains
(see $S$ and $T$ in (\ref{defST}))~:
\begin{equation}
  \label{eq:2023}
  \left( \vphantom{\frac12}
T(x,p,q,{\lambda})+S(x,p,q,{\lambda})\dot{x}\right)(p_x-q)^2=0.
\end{equation}
The second factor must be zero because, 
if $T+S\dot{x}$ was identically zero as a function of
$u,\ldots,u^{(k-1)},x,v,\ldots$ $v^{(\ell-1)}$, then,
by Definition~\ref{defparam} (point \ref{def1}), \emph{all} solutions
$(x(t),y(t),z(t))$ of $(\ref{sys3})$ would satisfy
$T(x,y,z,\dot{y}-z\dot{x})+\dot{x}S(x,y,z,\dot{y}-z\dot{x})=0$ identically, and
this would imply that $S$ and $T$ are identically zero functions of 4
variables, but we supposed the contrary.
The relation $q=p_x$ implies
\begin{equation}
  \label{PPP}
  \lambda=\mathcal{P}=Fp+rp_{u^{(k-1)}}+v^{(\ell)}p_{v^{(\ell-1)}}+v^{(\ell+1)}p_{v^{(\ell)}}
\end{equation}
and (\ref{eq201}) then yields
$p_{xx}=g(x,p,p_x,\lambda)$, or, with $\gamma$ defined by (\ref{defgamma}),
\begin{equation}
  \label{eqq30}
  \lambda=\gamma(x,p,p_x,p_{xx})\ .
\end{equation}
Since neither $p$ nor $Fp$ nor $r$ depend on $v^{(\ell+1)}$, (\ref{PPP}) and (\ref{eqq30}) yield
$p_{v^{(\ell)}}=0$, \ie\  $p$ is a
function of $u,\ldots,u^{(k-1)},x$, $v,\ldots,v^{(\ell-1)}$ only.
Then (\ref{PPP}) and (\ref{eqq30}) imply (\ref{ducon0}) with
$f=\gamma$. Furthermore, since $\varphi_{v^{(\ell)}}\neq0$ (point~\ref{def3} of Definition~\ref{defparam}),
(\ref{defr}) implies $r_{v^{(\ell)}}\neq0$.
Also, if $p$ was a function of $x$ only, then all solutions of (\ref{sysST0}) should
satisfy a relation $y(t)=p(x(t))$, which is absurd from Lemma~\ref{lem-absurd}.
We may then apply lemma~\ref{lem0} and assert that $k\geq1$, $\ell\geq1$,
$p_{u^{(k-1)}}\neq0$, $p_{v^{(\ell-1)}}\neq0$.

Since $p$ does not depend on $v^{(\ell)}$, (\ref{eqq30}) implies that the
right-hand side of (\ref{PPP}) does not depend on $v^{(\ell)}$ either; since
$p_{u^{(k-1)}}\neq0$, $r$ must be affine with respect to $v^{(\ell)}$, \ie
\begin{equation}
  \label{r}
  r\ =\ \tau\,+\,\sigma\,v^{(\ell)}\ ,
\end{equation}
with $\sigma$ and $\tau$ some functions of
$u,\ldots,u^{(k-1)},x,v,\ldots,v^{(\ell-1)}$.
Since $p$, $q=p_x$, $\lambda$ and $q_x=p_{xx}$ do not depend on $v^{(\ell)}$,
(\ref{eq2}) implies that ${\mathcal{Q}}$ does not depend on $v^{(\ell)}$
either; with $p_x=q$, and $r$ given by (\ref{r}), the expression of 
${\mathcal{Q}}_{v^{(\ell)}}$ is
$\sigma p_{xu^{(k-1)}}+p_{xv^{(\ell-1)}}$ while, from (\ref{PPP}), the expression of 
${\mathcal{P}}_{v^{(\ell)}}$.
Collecting this, one gets
\begin{equation}
  \label{eqq31}
  \sigma p_{u^{(k-1)}}+p_{v^{(\ell-1)}}=\sigma
p_{xu^{(k-1)}}+p_{xv^{(\ell-1)}}=0\ .
\end{equation}

Since $p_{u^{(k-1)}}\neq0$ and $p_{v^{(\ell-1)}}\neq0$,
(\ref{eqq31}) implies $E\sigma=0$, and also $\sigma_x=0$, $\sigma\neq0$. Then,
since $r_x\neq0$ (see (\ref{defr})), (\ref{r}) implies $\tau_x\neq0$.
With the above remarks, (\ref{PPP}) yields $\mathcal{P}=\lambda=Fp+\tau p_{u^{(k-1)}}$ 
and hence, from (\ref{eqq30}), the first relation in (\ref{EDPbis}).
In a similar way, (\ref{PQ}) yields $\mathcal{Q}=Fp_x+\tau p_{xu^{(k-1)}}$,
and substituting in (\ref{eq2}), one obtains (the terms involving $\dot{x}$
disappear according to (\ref{eqq30}))
$Fp_x-\delta(x,p,p_x,p_{xx})+p_{xu^{(k-1)}}\tau=0$ with $\delta$ defined by
(\ref{defdelta}).
This proves that $p$ satisfies (\ref{EDPbis}), equivalent to (\ref{EDPp}) according to
Remark~\ref{rmk-EDPbis}, and hence that $p$ is a solution of $\sysedp$.

To prove by contradiction that it is $K$-regular for some $K\leq k+\ell+1$,
assume that $ED^ip=0$ for $0\leq i\leq k+\ell$.
Then 
$p_x,p,\ldots,D^{(k+\ell-1)}p, x,\ldots,x^{(k+\ell-1)}$ are $2k+2\ell+1$
functions in the $2k+2\ell$ variables $u,\ldots,u^{(k-1)}, v,\ldots, v^{(\ell-1)},
x,\ldots,x^{(k+\ell-1)}$. At points where the Jacobian matrix has constant
rank, there is at least one nontrivial relation between them. From
point~\ref{def1} of Definition~\ref{defparam}, this would imply that all
solutions of system (\ref{sys3}) satisfy this relation, say
$R(z(t),y(t),\ldots,y^{(k+\ell-1)}(t),x(t),\ldots,x^{(k+\ell-1)}(t))=0$, which
is absurd from Lemma~\ref{lem-absurd}.
\end{proof}

\subsection{The case where $S$ and $T$ are zero}

Here, the situation is slightly more complicated:
we also establish that any parameterization ``derives from'' a solution of
the system of PDEs (\ref{EDPp}), but this is correct only if $J$ is not zero,
and there are two distinct (non equivalent) choices for $\gamma$ and $\delta$.
If $J\neq0$, we saw, in
section~\ref{sec-ST0}, that possibly after a change of coordinates, system
(\ref{sys3}) can be written as (\ref{eq:bof}), which we re-write here without
the tildes:
\begin{equation}
  \label{sysST0}
  \dot{z}=\kappa(x,y,z)\left(\dot{y}-\alpha(x,y,z)\,\dot{x}\right)
\left(\dot{y}-\beta(x,y,z)\,\dot{x}\right)
\;+\;a(x,y,z)\,\dot{x}+b(x,y,z)\,\dot{y}+c(x,y,z)\;,
\end{equation}
where $\kappa,\alpha,\beta,a,b,c$ are real analytic functions of three
variables and
$
\kappa\neq0$, $\alpha-\beta\neq0$, $\partial\alpha/\partial x\neq0$, $\partial\beta/\partial x\neq0$.
We state the theorem for this class of systems, because it is simpler to
describe the two possible choices for $\gamma$ and $\delta$ than with
(\ref{sys3}), knowing that $S=T=0$.

\begin{lmm}
If system (\ref{sysST0}) admits a parameterization $(\varphi,\psi,\chi)$ of order $(k,\ell)$ 
at a point, then $\varphi_{u^{(k)}}$ is a nonzero real analytic function.
\end{lmm}
\begin{proof}After a change of coordinates (\ref{diffST0}), use Lemma~\ref{lemphi}.\end{proof}

\begin{thrm}\label{th0}
Let $(x_0,y_0,z_0)$ be a point where $\kappa$, $\alpha-\beta$,
$\alpha_3$ and $\beta_3$ are nonzero, and $k,\ell,L$ three integers. 
If system (\ref{sysST0}) has a parameterization of order
$(k,\ell)$ at $\mathcal{X}=(x_0,y_0,z_0,\dot{x}_0,\dot{y}_0$, $\ldots,x^{(L)}_0$,
$y^{(L)}_0)$
with $\varphi,\psi,\chi$ defined on $U\subset\RR^{k+\ell+2}$,
then $k\geq1$, $\ell\geq 1$ and, for any point
$(u_0,\ldots,u_0^{(k)},v_0$, $\ldots,v_0^{(\ell)})\in U$ (not necessarily sent to
$\mathcal{X}$ by the parameterization) such that
$$\varphi_{u^{(k)}}(u_0,\ldots,u_0^{(k)},v_0,\ldots,v_0^{(\ell)})\neq 0,$$ 
there exist a neighborhood $\ouvedp$ of $(u_0,\ldots,u^{(k-1)}_0,\,\varphi(u_0\,\cdots\,v_0^{(\ell)})\,,v_0,\ldots,v^{(\ell-1)}_0)$
in $\RR^{k+\ell+1}$ and a regular solution  $p:\ouvedp\to\RR$ of one of the two systems
  $\sysedpl{\gamma^1}{\delta^1}k\ell$ or $\sysedpl{\gamma^2}{\delta^2}k\ell$ 
  with $\gamma^i$, $\delta^i$ given by (\ref{case1})-(\ref{case2}), such that $p,\varphi,\psi,\chi$
are related  by (\ref{defphi}), (\ref{defpsi}) and (\ref{defchi}).
\end{thrm}

Remark~\ref{rmk-Lmin} applies to this theorem in the same way as theorem~\ref{thSTnot0}.

\begin{proof}
Like in the beginning of the proof of Theorem~\ref{thSTnot0},
a parameterization $(\varphi, \psi, \chi)$
of order $(k,\ell)$ with $\varphi_{u^{(k)}}\neq0$ yields
an implicit form (\ref{eqp}).
Substituting in (\ref{sysST0}), one obtains an identity between two polynomials in
$v^{(\ell+1)}$ and $\dot{x}$.
The coefficient of $(v^{(\ell+1)})^2$ in the right-hand side must be zero and
this yields that $p$ cannot depend on $v^{(\ell)}$; the linear term in
$v^{(\ell+1)}$ then implies that $q$ does not depend on $v^{(\ell)}$ either.
To go further, let us define, as in the proof of Theorem~\ref{thSTnot0},
\begin{equation}
\label{j}
  \mathcal{P}=Fp+rp_{u^{(k-1)}}+v^{(\ell)}p_{v^{(\ell-1)}}\,,\ \ 
  \mathcal{Q}=Fq+rq_{u^{(k-1)}}+v^{(\ell)}q_{v^{(\ell-1)}}\,,
\end{equation}
with $F$ as in (\ref{F}). Still substituting in (\ref{sysST0}), the terms of degree 0, 1 and 2 with respect
to $\dot{x}$ then yield
\begin{equation} \label{eq:deg2xter}
\begin{array}{l}
{\mathcal{Q}}=\kappa(x,p,q){\mathcal{P}}^2+b(x,p,q){\mathcal{P}}+c(x,p,q)\,,\\
q_x=\kappa(x,p,q)\left(2p_x-\alpha(x,p,q)-\beta(x,p,q)\right)\mathcal{P}+a(x,p,q)+b(x,p,q)p_x\,,\\
0=\left(p_x-\alpha(x,p,q)\right)\left(p_x-\beta(x,p,q)\right)\,.
\end{array}\end{equation}
The factors in 
the third equation cannot both be zero because $\alpha-\beta\neq0$. Let us assume
\begin{equation}
  \label{pxqALT}
  p_x-\alpha(x,p,q)=0\;,\ \ \ 
p_x-\beta(x,p,q)\neq 0
\end{equation}
(interchange the roles of $\alpha$ and $\beta$ for the other alternative).
Since $\alpha_3\neq0$, the map $A$ defined in (\ref{diffST0}) has locally an
inverse $A^{-1}$, and the equation in (\ref{pxqALT}) is equivalent to
$(x,p,q)=A^{-1}(x,p,p_x)$; by differentiation an expression of $q_x$ as a
function of $x,p,p_x,p_{xx}$ is obtained; solving the second equation in 
(\ref{eq:deg2xter}) for $\mathcal{P}$ and substituting $q$ and $q_x$, 
one obtains $\mathcal{P}=\gamma^1(x,p,p_x,p_{xx})$ with $\gamma^1$ defined by
(\ref{case1})-(\ref{case2}). 
If one had chosen the other alternative in (\ref{pxqALT}), $A$ and $\gamma^1$
would be replaced by $B$ and $\gamma^2$.

Since $\mathcal{P}$ is also given by (\ref{j}), the relation (\ref{ducon0})
holds with $f=\gamma^1$;
also, for the same reasons as in the proof of Theorem~\ref{thSTnot0} (two
lines further than (\ref{eqq30})), $r_{v^{(\ell)}}$ is nonzero and it would be absurd
that $p$ depends on $x$ only.
One may then apply Lemma~\ref{lem0} and deduce that $k\geq1$, $\ell\geq1$,
$p_{u^{(k-1)}}\neq0$, $p_{v^{(\ell-1)}}\neq0$.

Since neither $p$ nor $\mathcal{P}=\gamma^1(x,p,p_x,p_{xx})$ depend on $v^{(\ell)}$
and $p_{u^{(k-1)}}\neq0$, the first equation in (\ref{j}) implies that $r$
assumes the form (\ref{r}) with $\sigma$ and $\tau$ some functions of the
$k+\ell+1$ variables
$u,\dot{u},\ldots,u^{(k-1)},x,v,\dot{v},\ldots,v^{({\ell-1})}$, and that two
relations hold: on the one hand $\sigma p_{u^{(k-1)}}+p_{v^{(\ell-1)}}=0$, i.e. one of the
relations in (\ref{EDPbis}), and on the other hand the first relation in
(\ref{EDPbis}) with $\gamma=\gamma^1$.
Similarly, the second equation in (\ref{j}) yields 
$\sigma q_{u^{(k-1)}}+q_{v^{(\ell-1)}}=0$ and 
$Fq+\tau q_{u^{(k-1)}}=\mathcal{Q}=\kappa\mathcal{P}^2+b\mathcal{P}+c$. 
Applying $F+\tau\partial/\partial{u^{(k-1)}}$ and $E$ to the first relation in
(\ref{pxqALT}) and using the four relations we just established, one obtains
on the one hand the second relation in (\ref{EDPbis}), with $\delta=\delta^1$
($\delta^1$ defined in (\ref{case1})-(\ref{case2})) and on the other hand 
$\sigma p_{xu^{(k-1)}}+p_{xv^{(\ell-1)}}=0$. The relations $\sigma_x=0$,
$\sigma\neq0$ and $\tau_x\neq0$ are then obtained exactly like at the end of
the proof of theorem~\ref{thSTnot0}; hence $p$ satisfies (\ref{EDPbis}) with
$\gamma=\gamma^1$ and $\delta=\delta^1$; this proves, thanks to
Remark~\ref{rmk-EDPbis}, that $p$ is a solution of
$\sysedpl{\gamma^1}{\delta^1}k\ell$ (it would be $\sysedpl{\gamma^2}{\delta^2}k\ell$
if one had chosen the other alternative in (\ref{pxqALT})).
The last paragraph of the proof of Theorem~\ref{thSTnot0} can be used to prove
that this solution is $K$-regular with $K\leq k+\ell+1$.
\end{proof}

\section{Flat outputs and differential flatness}
\label{sec-flat}

\begin{dfntn}[flatness, endogenous parameterization~\cite{Flie-Lev-Mar-R92cras}]
\label{def-flat}
  A pair $A\!=\!(a,b)$ of real analytic functions
  on a neighborhood of $(x_0,y_0,z_0,\ldots,x^{(j)}_0,y^{(j)}_0)$ in
  $\widehat{\Omega}\times\RR^{2\ordreplat-2}$ is a
  \emph{flat output of order $\ordreplat$} at
  $\mathcal{X}=(x_0,y_0,z_0,\ldots,x^{(L)}_0,y^{(L)}_0)$
  (with $L\geq j\geq0$) for system (\ref{sys3}) if there exists a Monge
  parameterization (\ref{para1}) of some order $(k,\ell)$ at $\mathcal{X}$
  such that any germ
  $(x(.),y(.),z(.),u(.),v(.))\in\mathcal{V}\times\mathcal{U}$ 
  (with $U,V$ possibly smaller than in (\ref{para1}))
  satisfies (\ref{flat1}) if and only if it satisfies (\ref{flat2}):
  \begin{eqnarray}
    \label{flat1}&\left.
    \begin{array}{rcl}
      \varphi\bigl(u(t),\dot{u}(t),\ldots,u^{(k)}(t),v(t),\dot{v}(t),\ldots,v^{(\ell)}(t)\bigr)&=&x(t)
      \\ 
      \psi\bigl(u(t),\dot{u}(t),\ldots,u^{(k)}(t),v(t),\dot{v}(t),\ldots,v^{(\ell)}(t)\bigr)&=&y(t)
      \\ 
      \chi\bigl(u(t),\dot{u}(t),\ldots,u^{(k)}(t),v(t),\dot{v}(t),\ldots,v^{(\ell)}(t)\bigr)&=&z(t)
    \end{array}
\right\}\;,
\\
    \label{flat2}&
\left.\!\!\!\!\!\!\!\!\!\!\!\!
    \begin{array}{l}
\dot{z}(t)=h\bigl(x(t),y(t),z(t),\,\dot{y}(t)\!-\!z(t)\dot{x}(t)\,\bigr)
\;+\;g\bigl(x(t),y(t),z(t),\,\dot{y}(t)\!-\!z(t)\dot{x}(t)\,\bigr)\;\dot{x}(t) 
      \\
      u(t)=a\bigl(x(t),y(t),z(t),\dot{x}(t),\dot{y}(t),\ddot{x}(t),\ddot{y}(t),
      \ldots, x^{(\ordreplat)}(t),y^{(\ordreplat)}(t)\bigr)
      \\ 
      v(t)=b\bigl((x(t),y(t),z(t),\dot{x}(t),\dot{y}(t),\ddot{x}(t),\ddot{y}(t),
      \ldots, x^{(\ordreplat)}(t),y^{(\ordreplat)}(t)\bigr)\
    \end{array}
\right\}\!.
  \end{eqnarray}
System (\ref{sys3}) is called \emph{flat} if and only if it admits a flat
output of order $\ordreplat$ for some $\ordreplat\in\NN$.
A Monge parameterization is \emph{endogenous}\footnote{
This terminology (endogenous vs. exogenous) is borrowed from the authors of
\cite{Flie-Lev-Mar-R92cras,Mart92th}; it usually qualifies feedbacks rather than
  parameterizations, but the notion is exactly the same.
}
if and only if there exists a flat output associated to this parameterization
as above.
\end{dfntn}
In control theory, flatness is a better
known notion than Monge parameterization.
For general control systems, it implies existence of a parameterization
(obvious in the above definition), and people 
conjecture~\cite{Flie-Lev-Mar-R99open} that the two notions are in fact
equivalent, at least away from some singular points. In any case, our results are relevant to both: systems (\ref{sys3}) that are
proved to be parameterizable are also flat and our efforts toward proving that
the other ones are not parameterizable would also prove that they are not flat.

Theorem~\ref{edpparam} gave a procedure to derive a parameterization of
(\ref{sys3}) from a regular solution $p$ of $\sysedp$, and we saw in
Section~\ref{sec-main} that, unless $S=T=J=0$, these are the only possible
parameterizations. One can tell when such a parameterization is endogenous:

\begin{prpstn}
  Let $p:\ouvedp\to\RR$, with $\ouvedp\subset\RR^{k+\ell+1}$ open,  be a
regular solution of system $\sysedp$.
The parameterization of order $(k,\ell)$ 
of system (\ref{sys3}) associated to $p$ according to Theorem~\ref{edpparam}
is endogenous if and only if $p$ is exactly $(k+\ell-2)$-regular; then, the
associated flat output is of order $\ordreplat\leq k+\ell-2$.
\end{prpstn}
\begin{proof}
  In the end of the proof of Theorem~\ref{edpparam}, it was established that
  (\ref{flat1}), written 
  $\Gamma(u,v)=(x,y,z)$, is equivalent to (\ref{uvij}) if $(x,y,z)$ is a
  solution of (\ref{sys3}). 
If either $i_0<k$ or $j_0<\ell$ in (\ref{uvij}), then there are, for fixed
  $x(.),y(.),z(.))$, infinitely many solutions $u(.),v(.))$ of (\ref{uvij})
  while there is a unique one for (\ref{flat2}). Hence $i_0=k$ and $j_0=\ell$
  if (\ref{flat1}) is equivalent to (\ref{flat2}); then $K=i_0+j_0-2=k+\ell-2$
  so that $p$ is $(k+\ell-2)$-regular and (\ref{uvij}) (where $u$ and $v$ do
  not appear in the right-hand side) is of the form (\ref{flat2}) with
  $\ordreplat=K=k+\ell-2$. 
\end{proof}
The main result in~\cite{Pome97cocv} is a necessary condition for ``$(x,u)$-dynamic linearizability'' ($(x,u)$-flatness
might be more appropriate) of system (\ref{sys4}).
For system (\ref{sys4}), it means existence of a flat output whose components
are functions 
of $\xi^1,\xi^2,\xi^3,\xi^4,w^1,w^2$; for system (\ref{sys3}), it translates
as follows.
The functions $\gamma$ and $\delta$ in
(\ref{sys4}) are supposed to be related to  $g$ and $h$ in
(\ref{sys3}) according to (\ref{defgamma}) and (\ref{defdelta}).
\begin{dfntn}
\label{def-xu}
System (\ref{sys4}) is ``$(x,u)$-dynamic linearizable'' is and only if system (\ref{sys3}) admits a flat output of order 2 of a
special kind~: 
$A(x,y,z,\dot{x},\dot{y},\ddot{x},\ddot{y})=\mathfrak{a}(x,y,z,\lambda,\dot{x},\dot{\lambda})$
for some smooth $\mathfrak{a}$. 
\end{dfntn}
The following proposition is usefull to recover the main result from
\cite{Pome97cocv}, see Remark~\ref{rmrk-old}. 
\begin{prpstn}
\label{prop-xu33}
If   system \textup{(\ref{sys4})} is ``$(x,u)$-dynamic linearizable'' in the sense of
  \textup{\cite{Pome97cocv}}, then \textup{(\ref{sys3})} admits a parameterization of order
  $(k,\ell)$ with $k\leq3$ and $\ell\leq3$.
\end{prpstn}
\begin{proof}$\ $\\[-1\baselineskip]
Consider the map
$\displaystyle
  (x,y,z,\lambda,\dot{x},\dot{\lambda},\ldots,x^{(4)},\lambda^{(4)}) \mapsto
  \left( \!\!
    \begin{array}{c}
\mathfrak{a}(x,y,z,\lambda,\dot{x},\dot{\lambda})
\\
\dot{\mathfrak{a}}(x,y,z,\lambda,\dot{x},\dot{\lambda},\ddot{x},\ddot{\lambda})
\\
\ddot{\mathfrak{a}}(x,y,z,\lambda,\dot{x},\dot{\lambda},\ldots,x^{(3)},\lambda^{(3)})
\\
\mathfrak{a}^{(3)}(x,y,z,\lambda,\dot{x},\dot{\lambda},\ldots,x^{(4)},\lambda^{(4)})
    \end{array}
  \!\!\right).
$ 
\\
Its Jacobian is $8\times12$, and has rank 8, but the $8\times8$ sub-matrix
corresponding to derivatives with respect to
$\dot{x},\dot{\lambda},\ldots,x^{(4)},\lambda^{(4)}$ has rank 4 only.
Hence $x$, $y$, $z$, and $\lambda$ can be expressed as functions of the
components of
$\mathfrak{a},\dot{\mathfrak{a}},\ddot{\mathfrak{a}},\mathfrak{a}^{(3)}$,
yielding a Monge parameterization of order at most $(3,3)$.
  \end{proof}

\section{Conclusion}
\label{sec-concl}

Let us discuss both flatness (see Section~\ref{sec-flat}) and Monge parameterization.
For convenience, assume $k\leq\ell$ 
and call F-systems the systems (\ref{sys3}) such
that $S=T=J=0$ and C-systems all the other ones.

\medskip

F-systems are flat; this was proved in~\cite{Pome97cocv}. This paper adds that they 
admit a Monge parameterization of order (1,2), but does \emph{not} prove 
differential flatness of any system not known to be flat up to now:
C-systems are not believed to be flat.
It does not either prove non-flatness of any system: 
it only \emph{conjectures} that no C-system admits a parameterization, and hence none of
them is flat. To the best of our
knowledge, no one knows whether simple systems like
(\ref{eq:ex}-b) or (\ref{eq:ex}-c) are flat of not.

The first contribution of the paper is to prove that a C-system admits a
parameterization of order $(k,\ell)$ if and only if the PDEs
$\sysedp$, for suitable $\gamma,\delta$, admit a ``regular solution'' $p$. The second
contribution is to prove that, for any $\gamma,\delta$, there is no regular
solution to $\sysedp$ if either $k\leq2$ or $k=\ell=3$
(this does not contradict existence of parameterizations of order $(1,2)$ for
F-systems: these do not ``derive from'' a solution of these PDEs). We guess, in
Conjecture~\ref{conj-EDP}, that even for higher values of the integers $k,\ell$,
\emph{none} of these PDEs have any regular solution; this would
imply that C-systems are not flat.

Besides recovering the results from~\cite{Pome97cocv} with far more
natural and elementary arguments, we believe that some insight was gained on
Monge parameterizations of 
\emph{any order} for ``C-systems'',  by reducing non-parameterizability to
non-existence of solutions to a systems of PDEs that can easily be written for
any $k,\ell$.

\medskip

The main perspective raised by this paper is to \textbf{prove Conjecture~\ref{conj-EDP}}.
The only \emph{theoretical} difficulty is, in fact, that no a priori bound on the
integers $k,\ell$ is known. 
Indeed, as explained in Section~\ref{sec-EDP}, for fixed
$k,\ell,\gamma,\delta$, it amounts to a classical problem.
To prove Proposition~\ref{prop-edp}, we solved, in a synthetic manner, that
problem for $k\leq2$ or $k=\ell=3$ and arbitrary $\gamma$ and $\delta$.
We lack a non-finite argument, or a better understanding of the
structure, to go to arbitrary $k,\ell$. 
Let us comment more on the (non trivial) case where
$\gamma$ and $\delta$ are \emph{polynomials}, for instance the very simple ones in
(\ref{EDPp2}). For fixed $k,\ell$, the question can be formulated in terms of
differential polynomial rings:
does the differential ideal generated by left-hand sides of the equations (\ref{EDPp2})
contain the polynomials $ED^ip$~?
Differential elimination (see~\cite{Ritt50} or the recent
survey~\cite{Hube03}) is relevant here; finite algorithms
have been already implemented in computer algebra.
Although we have not yet succeeded (because of complexity) in carrying out
these computations, even on example (\ref{EDPp2}) for $(k,\ell)=(3,4)$, and
although it will certainly not \emph{provide} a bound on $k,\ell$, we do believe that 
computer algebra is a considerable potential help.

\smallskip
Another perspective is to enlarge the present approach to \textbf{higher
  dimensional control systems}. 
For instance, what would play the role of our system of PDEs $\sysedp$ when,
instead of (\ref{sys3}), one considers a single relation between more than
three scalar functions of time (this captures, instead of (\ref{sys4}),
control affine systems with $n$ states and 2 controls, $n>4$)~? 
We have very little insight on this question: the present paper strongly takes
advantage of the special structure inherent to our small dimension; the
situation could be far more complex.


\clearpage
\appendix
\section{}
\label{app-1}

\begin{proof}[Proof of Lemma \ref{jacob}]
For this proof only, the notation $\mathcal{F}_{i,j}$ ($0\leq i\leq k$, 
$0\leq j\leq\ell$) stands either for the following family of $i+j$ vectors in
$\RR^{K+2}$ or for the corresponding $(K+2)\times(i+j)$ matrix~:
\begin{displaymath}
  \mathcal{F}_{i,j}\ =\ \left(
\partialxy{\pi}{u^{(k-i)}},\ldots,\partialxy{\pi}{u^{(k-1)}},
\partialxy{\pi}{v^{(\ell-j)}},\ldots,\partialxy{\pi}{v^{(\ell-1)}}
\right)
\end{displaymath}
with the convention that if $i$ or $j$ is zero the corresponding list is empty;
$\mathcal{F}_{i,j}$ depends on
$u,\ldots,u^{(k-1)},v,\ldots$, $v^{(\ell-1)},x,\ldots,x^{(K)}$.
Let us first prove that, at least outside a closed subset of empty interior,
\begin{equation}
  \label{eq:hh}
  \rank\mathcal{F}_{k,\ell}\ =\ K+2\ .
\end{equation}
Indeed, if it is smaller at all points of $\ouvedp\times\RR^K$, then, around
points (they form an open dense set) where it is locally constant, there
is at least one function $R$ such
that a non-trivial identity
$R(p_x,p,\ldots,D^{K}p,x,\ldots,x^{(K)})=0$ holds and the partial
derivative of $R$ with respect to at least one of its $K+2$ first arguments is
nonzero.
Since $p$ is $K$-regular, applying $E$ to this relation, shows that $R$
does not depend on $D^{K}p$, and hence does not depend on $x^{(K)}$
either. Then, applying $ED$, $ED^2$ and so on, and using the fact that,
according to (\ref{psolgh}), $Dp_x$ is a function of $p_x,p,Dp,x,\dot{x}$, we get finally a relation $R(p_x,p,x)=0$ with
$(R_{p_x},R_p)\neq(0,0)$.
Differentiating with respect to $u^{(k-1)}$, one obtains
$R_{p_x}p_{xu^{(k-1)}}+R_p p_{u^{(k-1)}}=0$; hence, from the first relation in
(\ref{EDPp}-c), $R_{p_x}\neq 0$, and the relation $R(p_x,p,x)=0$ implies, in a
neighborhood of almost any point, $p_x=f(p,x)$
for some smooth function $f$. From Lemma~\ref{lemAAA}, this would contradict
the fact that the solution $p$ is $K$-regular. This proves (\ref{eq:hh}).

\smallskip
Let now $W_s$ ($1\leq s\leq K+2$) be the set of pairs $(i,j)$ such that $i+j=s$
and the rank of $\mathcal{F}_{i,j}$ is $s$ at least at one point in
$\ouvedp\times\RR^K$, \ie\  one of the $s\times s$ minors of $\mathcal{F}_{i,j}$ is a
nonzero real analytic function on $\ouvedp\times\RR^K$. The lemma states that $W_{K+2}$ is
nonempty; in order to prove it by contradiction, suppose that
$W_{K+2}=\varnothing$ and
let $\overline{s}$ be the smallest $s$ such that $W_s=\varnothing$. 
From (\ref{EDPp}-c), $W_1$ contains $(1,0)$, hence $2\leq \overline{s}\leq K+2<k+\ell+1$. 
Take $(i',j')$ in $W_{\overline{s}-1}$; $\mathcal{F}_{i',j'}$ has rank $i'+j'$ (\ie\  is
made of $i'+j'$ linearly independent vectors) on an open dense set 
$A\subset \ouvedp\times\RR^K$.
Let the $i_1\leq k$ and $j_1\leq\ell$ be the largest such that
$\mathcal{F}_{i_1,j'}$ and $\mathcal{F}_{i',j_1}$ have rank $\overline{s}-1$ on $A$. 
On the one hand, since $i'+j'=\overline{s}-1<k+\ell$, one has either
$i'<k$ or $j'<\ell$. On the other hand since $W_{\overline{s}}$ is empty, it contains neither $(i'+1,j')$ nor
$(i',j'+1)$; hence the rank of $\mathcal{F}_{i'+1,j'}$ is less than $i'+j'+1$
if $i'<k$, and so is the rank of $\mathcal{F}_{i',j'+1}$ if $j'<\ell$. 

To sum up, the following implications hold:
$\;i'<k\Rightarrow i_1\geq i'+1\;$ and $\;j'<\ell\Rightarrow j_1\geq j'+1\;$.
From (\ref{eq:hh}), one has either $i_1<k$ or $j_1<\ell$.
Possibly exchanging $u$ and $v$, assume  $i_1<k$; all the vectors
$\partial{\pi}/\partial{u^{(k-i_1)}}$, \ldots,
  $\partial{\pi}/\partial{u^{(k-i'+1)}}$, 
$\partial{\pi}/\partial{u^{(\ell-j_1)}}$, \ldots,
$\partial{\pi}/\partial{u^{(\ell-j'+1)}}$ are then linear combinations of the vectors in
$\mathcal{F}_{i',j'}$, while $\partial{\pi}/\partial{u^{(k-i_1-1)}}$ is not~:
\begin{equation}
  \label{eqq1}
  \rank\mathcal{F}_{i',j'}=i'+j'\,,\ \ \ \ 
  \rank\mathcal{F}_{i_1,j_1}=i'+j'\,,\ \ \ \ 
\rank\left(
\partialxy{\pi}{u^{(k-i_1-1)}}\,,\,\mathcal{F}_{i',j'}
\right)=i'+j'+1
\end{equation}
on an open dense subset of $\ouvedp\times\RR^K$, that we still call $A$ although it
could be smaller. 
In a neighborhood of any point in this set, one can, from the third relation, apply
the inverse function theorem and obtain, for an open
$\Omega\subset\RR^{k+\ell+K+1}$, a map $\Omega\to\RR^{i'+j'+1}$ that expresses
$u^{(k-i')},\ldots,u^{(k-1)}$, $v^{(\ell-j')},\ldots,v^{(\ell-1)}$ and
$u^{(k-i_1-1)}$ as functions of
$u,\ldots,u^{(k-i_1-2)}$, $u^{(k-i_1)},\ldots,u^{(k-i'-1)}$, $v,\ldots,v^{(\ell-j'-1)}$,
$x,\ldots,x^{(K)}$ and $i'+j'+1$ functions
chosen among $p_x,p,Dp,\ldots,D^{K}p$ ($i'+j'+1$
columns defining an invertible minor in $\left(
\partial{\pi}/\partial{u^{(k-i_1-1)}}\,,\,\mathcal{F}_{i',j'}
\right)$). 
Focusing on $u^{(k-i_1-1)}$, one has
\begin{equation}
\label{eqq5}
\!\!\!\!\!\!\begin{array}{l}
\!\!u^{(k-i_1-1)}=\\
B\left(u,\ldots,u^{(k-i_1-2)},u^{(k-i_1)},\ldots,u^{(k-i'-1)},v,\ldots,v^{(\ell-j'-1)},x,\ldots,x^{(K)},p_x,p,\ldots,D^{K}p\right)
\end{array}
\end{equation}
where $B$ is some smooth function of $k+\ell+2K+2-i'-j'$ variables and
we have written all the functions $p_x,p,Dp,\ldots,D^{K-1}p$ although $B$
really depends only on $i'+j'+1$ of them.

Differentiating (\ref{eqq5}) with respect to
$u^{(k-i')},\ldots,u^{(k-1)},v^{(\ell-j')},\ldots,v^{(\ell-1)}$, one has, with
obvious matrix notation, $\left(\partialxy{B}{p_x}\ \partialxy{B}{p}\ \cdots\
  \partialxy{B}{D^{K-1}p}\right)\mathcal{F}_{i',j'}=0$,
where the right-hand side 
is a line-vector of dimension $i'+j'$; from (\ref{eqq1}), this implies 
\begin{equation}
  \label{eqq60}
  \left(\partialxy{B}{p_x}\ \partialxy{B}{p}\ \cdots\
  \partialxy{B}{D^{K-1}p}\right)\mathcal{F}_{i_1,j_1}\ =\ 0\;,
\end{equation}
where the right-hand side 
is a now a bigger line-vector of dimension $i_1+j_1$.
Differentiating (\ref{eqq5}) with respect to
$u^{(k-i_1)},\ldots,u^{(k-i'-1)},v^{(\ell-j_1)},\ldots,v^{(\ell-j'-1)}$ and
using (\ref{eqq60}) yields that $B$ does not depend on its arguments
$u^{(k-i_1)},\ldots,u^{(k-i'-1)}$ and $v^{(\ell-j_1)},\ldots,v^{(\ell-j'-1)}$. 
$B$ cannot depend on $D^{K}p$ either because $ED^{K}p\neq0$ and all the
other arguments of $B$ are constant along $E$; then it cannot depend on $x^{(K)}$ either because
$x^{(K)}$ appears in no other argument; (\ref{eqq5}) becomes
\begin{displaymath}
  u^{(k-i_1-1)}=B\left(u,\ldots,u^{(k-i_1-2)},v,\ldots,v^{(\ell-j_1-1)},x,\ldots,x^{(K-1)},p_x,p,\ldots,D^{K-1}p\right)\;.
\end{displaymath}
Applying $D$, using (\ref{psolgh}) and
substituting $u^{(k-i_1-1)}$ from above, one gets, from some smooth $C$ ,
\begin{equation}
\label{eq:00001}
  u^{(k-i_1)}=C\bigl(\, u,\ldots,u^{(k-i_1-2)},
\underbrace{v,\ldots,v^{(\ell-j_1)}}_{\substack{{\rm empty\; if}\; j_1=\ell}},
x,\ldots,x^{(K)},p_x,p,\ldots,D^{K}p\bigr)\;.
\end{equation}
Differentiating with respect to
$u^{(k-i')},\ldots,u^{(k-1)},v^{(\ell-j')},\ldots,v^{(\ell-1)}$ yields\\
$\left(\partialxy{C}{p_x}\ \partialxy{C}{p}\ \cdots\
  \partialxy{C}{D^{K-1}p}\right)\mathcal{F}_{i',j'}=0$,
the right-hand side being a line-vector of dimension $i'+j'$. From the first
two relations in
(\ref{eqq1}), $\partial\pi/\partial u^{(k-i_1-1)}$ is a linear
combination of the columns of $\mathcal{F}_{i',j'}$, hence one also has
$\displaystyle 
\left(\partialxy{C}{p_x}\ \partialxy{C}{p}\ \cdots\
  \partialxy{C}{D^{K-1}p}\right)\partialxy{\pi}{u^{(k-i_1-1)}}\ =\ 0\ .
$ 
\\This implies that the derivative of
the right-hand side of (\ref{eq:00001}) with respect to $u^{(k-i_1-1)}$ is
zero. This is absurd. 
\end{proof}

\section{Proof of Lemmas \ref{lem-1}, \ref{lem-2} and \ref{lem-3}}
\label{app-sys}

We need some notations and
preliminaries.
With $F$, $E$ and $\tau$ defined in (\ref{F}) and (\ref{defsigtau}), define the vector fields
\begin{eqnarray}
  \label{XY}&\displaystyle
  \XX=\partialx{x}\;,\ \ \ 
\FF=F+\tau\partialx{u^{(k-1)}}\ 
\\
  \label{eq:bra}&
  \XX_1=[\XX,\FF],\ \XX_2=[\XX_1,\FF]\,,\hspace{2em}
\EE_2=[E,\FF],\ \EE_3=[\EE_2,\FF]\ .
\end{eqnarray}
Then (\ref{EDPbis}) obviously implies
\begin{equation}
  \label{YpYpx}
  Yp=\gamma(x,p,p_x,p_{x,x})\,,\ \ \ \ Yp_x=\delta(x,p,p_x,p_{x,x})\,,\ \ \ \ \XX\sigma=0\,,
\end{equation}
and a simple computation yields (we recall $E$ from (\ref{defsigtau}))~:
\begin{eqnarray}
\label{66}&
\XX_1=\tau_x\,\partialx{u^{(k-1)}}\,,\ \ 
\XX_2=\tau_x\,\partialx{u^{(k-2)}}+(\cdots)\partialx{u^{(k-1)}}\,,\ \ 
\\
\nonumber 
&
E=\partialx{v^{(\ell-1)}}+\sigma\,\partialx{u^{(k-1)}}\,,\ \ 
\EE_2=\partialx{v^{(\ell-2)}}+\sigma\,\partialx{u^{(k-2)}}
+(\cdots)\partialx{u^{(k-1)}}\,,\ \ 
\\
\nonumber 
&
\EE_3=\partialx{v^{(\ell-3)}}+\sigma\,\partialx{u^{(k-3)}}
+(\cdots)\partialx{u^{(k-2)}}+(\cdots)\partialx{u^{(k-1)}}\,.
\end{eqnarray}
The vector field $\XX_1$ and $\XX_2$ are linearly independent because
$\tau_x\neq0$, see (\ref{EDPbis}).
Computing the following brackets
and decomposing on $\XX_1$ and $\XX_2$, one gets
\begin{eqnarray}
  \label{bralie}&
  [\XX,\XX_1]=\lambda \XX_1\,,\ \ \   [\XX_1,\XX_2]=\lambda'\XX_1+\lambda''\XX_2\ ,
\\
  \label{bralie2}&
[\XX,\EE]=0\,,\ \ \ [\XX,\EE_2]=\mu\XX_1\,,\ \ \
[\XX,\EE_3]=\mu'\XX_1+\mu''\XX_2\,,\ \ \ 
[\EE_2,\XX_2]=\nu'\XX_1+\nu''\XX_2\ .
\end{eqnarray}
for some functions $\lambda,\lambda',\lambda'',\mu,\mu',\mu'',\nu',\nu''$.

\begin{proof}[Proof of Lemma \Rref{lem-1}]
From (\ref{EDPp}-c), $y=p(u,\ldots, u^{(k-1)},x,v,\ldots, v^{(\ell-1)})$ defines
local coordinates 
$u,\ldots u^{(k-2)}$, $y,x$, $v,\ldots,v^{(\ell-1)}$. 
Composing $p_x$ by the inverse of this change of coordinates, there is a function $\alpha$ of $k+\ell+1$ variables such that
$p_x=\alpha(u,\ldots,u^{(k-2)},\,p\,,x,\ldots, v^{(\ell-1)})$ identically. 
Since $Ep=Ep_x=0$ (see (\ref{EDPbis})), applying $E$ to both sides of this identity
yields that $\alpha$ does not depend on its argument $v^{(\ell-1)}$.
Similarly, if $k\geq2$, 
differentiating both sides of the same identity with respect to
$u^{(k-1)}$ and $u^{(k-2)}$, the fact that the determinant in the
lemma is zero implies that $\alpha$ does
not depend on its argument $u^{(k-2)}$.
To sum up, $p$ and $p_x$ satisfy an identity
$$p_x=\alpha(u,\ldots,u^{(k-3)},\,p\,,x,v,\ldots, v^{(\ell-2)}),$$
where the first list is empty if $k=1$ or $k=2$.
Now define two integers $m\leq k-3$ and $n\leq\ell-2$ as
the smallest such that
$\alpha$ depends on $u,\ldots,u^{(m)},x,y,v,\ldots,v^{(n)}$, with the
convention that $m<0$ if $k=1$, $k=2$, or $\alpha$ depends on none of the
variables $u,\ldots,u^{(k-3)}$ and $n<0$ 
if $\alpha$ depends on none of the variables $v,\ldots,v^{(\ell-2)}$.

Applying $\FF$ to both sides of the above identity yields
$$
\FF p_x=\alpha_y\,\FF p
+\sum_{i=0}^{m}u^{(i+1)}\alpha_{u^{(i)}}\,+\sum_{i=0}^{n}v^{(i+1)}\alpha_{
  v^{(i)}}\ ,
$$
where, if $m<0$ or $n<0$, the corresponding sum is empty.
Using (\ref{YpYpx}), since
$p_{xx}=\alpha_x+\alpha\alpha_y$, one can replace  $\FF p$ with
$\gamma(x,y,\alpha,\alpha_x+\alpha\alpha_y)$ and $\FF p_x$ with
$\delta(x,y,\alpha,\alpha_x+\alpha\alpha_y)$
in the above equation, where all terms except the last one of each non-empty
sum therefore depend on $u,\ldots,u^{(m)},x,y,v,\ldots,v^{(n)}$ only. Differentiating with
respect to $u^{(m+1)}$ and $v^{(n+1)}$ yields $\alpha_{u^{(m)}}=\alpha_{v^{(n)}}=0$,
which is possible only if $m<0$ and $n<0$, hence the lemma.
\end{proof}

\begin{proof}[Proof of Lemma \Rref{lem-2}]
From (\ref{inv2}), setting
\begin{equation}
  \label{yz}
y=p(u,\ldots,u^{(k-1)},x,v,\ldots,v^{(\ell-1)}),\ 
z=p_x(u,\ldots,u^{(k-1)},x,v,\ldots,v^{(\ell-1)}),
\end{equation}
one gets some local coordinates 
$(u,\ldots,u^{(k-3)},x,y,z,v,\ldots,v^{(\ell-1)})$.
In these coordinates, the vector fields
$\XX$ and $\FF$ defined by (\ref{XY}) have the following expressions, where
$\chi$ and $\alpha$ are some functions, to be studied further~:
\begin{eqnarray}
  \label{3-X}
  \XX&=&\partialx{x}+z\partialx{y}+\alpha\partialx{z}\;,
\\
  \label{3-F}
  \FF&=&\gamma\partialx{y}+\delta\partialx{z}
+\chi\partialx{u^{(k-3)}}
+\sum_{i=0}^{k-4}u^{(i+1)}\partialx{u^{(i)}}
+\sum_{i=0}^{\ell-1}v^{(i+1)}\partialx{v^{(i)}}\;.
\end{eqnarray}
In the expression of $\FF$, the third term is zero if $k=2$,
the fourth term ($\sum_{i=0}^{k-4}\cdots$) is zero if $k=2$ or $k=3$,
and the notations $\gamma$ and $\delta$ are slightly abusive~: $\gamma$
  stands for the function
$$(u,\ldots,u^{(k-3)},x,y,z,v,\ldots,v^{(\ell-1)})\mapsto\gamma(x,y,z,\alpha(u,\ldots,u^{(k-3)},x,y,z,v,\ldots,v^{(\ell-1)}))\;,$$
and the same for $\delta$. With the same abuse of notations, (\ref{EDPp}-e) reads
\begin{equation}
  \label{taux3}
  \XX\gamma-\delta\neq0.
\end{equation}
The equalities
$(\sigma\partialx{u^{(k-1)}}+\partialx{v^{(\ell-1)}})u^{(k-2)}
=\partialx{x}u^{(k-2)}=\partialx{u^{(k-1)}}u^{(k-2)} =0$ are obvious in the
original coordinates. 
Since the inverse of the change of coordinates (\ref{yz}) is given by
$$
u^{(k-2)}\!=\chi(u,\ldots,u^{(k-3)},x,y,z,v,\ldots,v^{(\ell-1)}),\ 
u^{(k-1)}\!=\FF\!\chi(u,\ldots,u^{(k-3)},x,y,z,v,\ldots,v^{(\ell-1)}),
$$
and  $\EE$, $\XX$ and $\XX_{1}$ are given by (\ref{66}), those  equalities imply
\begin{equation}
  \label{eq:666}
  \EE\chi=\XX\chi=\XX_{\,1}\,\chi=0\ .
\end{equation}
Then, from (\ref{eq:bra}), (\ref{3-X}) and (\ref{3-F}), 
\begin{eqnarray}
\label{cro1}
\XX_1&=& \left(\XX\gamma-\delta\right)\partialx{y} \;+\;
\left(\XX\delta-\FF\alpha\right)\partialx{z} \;,
\\[0pt]
  \label{cro2}
[\XX,\XX_1] &=& 
\left(\XX^2\gamma-2\XX\delta-\FF\alpha\right)\partialx{y} 
\;+\;
\left(\XX^2\delta-\XX\FF\alpha-\XX_1\alpha\right)\partialx{z}\ .
\end{eqnarray}
With these expressions of $\XX$ and $\XX_1$, the first relation in
(\ref{bralie}) implies~:
\begin{equation}
  \label{det0}
    \left|
    \begin{array}{cc}
\XX\gamma-\delta&\XX^2\gamma-2\XX\delta+\FF\alpha
\\
\XX\delta-\FF\alpha&\XX^2\delta-\XX\FF\alpha-\XX_1\alpha
    \end{array}
  \right|=0 \ .
\end{equation}

The definition of $\alpha$ implies $\XX z=\alpha$.
In the original coordinates, this translates into the identity
$p_{xx}=\alpha(u,\ldots,u^{(k-3)},x,p,p_x,v,\ldots,v^{(\ell-1)})$. 
Since $Ep=Ep_x=Ep_{x,x}=0$ (see (\ref{EDPbis})), applying $E$ to both sides of this identity
yields that $\alpha$ does not depend on its argument $v^{(\ell-1)}$. Also,
if $k\geq 3$, 
differentiating both sides with respect to $u^{(k-1)}$, $u^{(k-2)}$ and
$u^{(k-3)}$, we obtain that the determinant (\ref{eq:det2}) is zero if and
only if $\alpha$ does not depend on its argument $u^{(k-3)}$. 
To sum up, under the assumptions of the lemma,
\begin{equation}
  \label{al2}
  \alpha\ \mbox{depends on}\ \
  u,\ldots,u^{(k-4)},x,y,z,v,\ldots,v^{(\ell-2)}
\ \ \mbox{only}
\end{equation}
with the convention that the first list is empty if $k=2$ or $k=3$.
Now define two integers $m\leq k-4$ and $n\leq\ell-2$ as
the smallest such that
$\alpha$ depends on $u,\ldots,u^{(m)},x,y,v,\ldots,v^{(n)}$, with the
convention that $m<0$ if $k=2$, $k=3$, or $\alpha$ depends on none of the
variables $u,\ldots,u^{(k-4)}$, and $n<0$ 
if $\alpha$ depends on none of the variables $v,\ldots,v^{(\ell-2)}$.
We have
\begin{equation}
  \label{eq:mn}
  m\geq0\;\Rightarrow\;\alpha_{u^{(m)}}\neq0\ ,\ \ \ \ n\geq0\;\Rightarrow\;\alpha_{v^{(n)}}\neq0\ .
\end{equation}
Since $m$ is no larger that $k-4$, $\chi$ does not appear in 
the expression of $\FF\alpha$~:
\begin{eqnarray}
  \label{Falpha3}
  \FF\alpha & = & \gamma\alpha_{y}+\delta\alpha_{z}
+\sum_{i=0}^{m}u^{(i+1)}\alpha_{u^{(i)}}
+\sum_{i=0}^{n}v^{(i+1)}\alpha_{v^{(i)}}
\end{eqnarray}
where the first (or second) sum is empty if $m$ (or $n$) is negative.

In the left-hand side of (\ref{det0}), all the terms depend only on
$u,\ldots,u^{(m)},x,y,z,v,\ldots,v^{(n)}$, except $\FF\alpha$, $\XX\FF\alpha$
and $\XX_1\alpha$  that depend on
$u^{(m+1)}$ if $m\geq0$ or on $v^{(n+1)}$ if $n\geq0$ (see above); the
determinant is a polynomial of degree two with respect to
$u^{(m+1)}$ and $v^{(n+1)}$ with coefficients depending on $u,\ldots,u^{(m)},x,y,z,v,\ldots,v^{(n)}$ only, and the term of degree two, coming from $(\FF\alpha)^2$, is
$$
\left(\alpha_{u^{(m)}}u^{(m+1)}
+\alpha_{v^{(n)}}v^{(n+1)}\right)^2\ .
$$
Hence (\ref{det0}) implies $\alpha_{u^{(m)}}=\alpha_{v^{(n)}}=0$ and, from
(\ref{eq:mn}), negativity of $m$ and $n$ are negative. By definition of these integers, this implies that $\alpha$ depends on $(x,y,z)$
only: in the original coordinates, one has $p_{xx}=\alpha(x,p,p_x)$.
\end{proof}

Before proving Lemma~\ref{lem-3}, we need to extract more information from the
previous proof~:
\begin{lmm}
\label{lem-2cro}
Assume, as in Lemma~\ref{lem-2}, that $p$ is a solution of $\sysedp$ satisfying
(\ref{inv2}), but assume also that
$\ell\geq k\geq3$ and the determinant (\ref{eq:det2}) is nonzero.
Then $[\XX,\EE_2]=[\XX,\EE_3]=0$.
\end{lmm}
\begin{proof}
  Starting as in the proof of Lemma~\ref{lem-2}, one does not obtain
(\ref{al2}) but, since (\ref{eq:det2}) is nonzero,
  \begin{equation}
    \label{al2'}
  \alpha\ \mbox{depends on}\ \
  u,\ldots,u^{(k-4)},x,y,z,v,\ldots,v^{(\ell-2)}
\ \ \mbox{and}\ \     \alpha_{u^{(k-3)}}\neq 0\ .
  \end{equation}
Since $Ep=Ep_x=0$, one has $E=\partial/\partial v^{(\ell-1)}$ in these
coordinates.
The first equation in (\ref{eq:666}) then
reads $\chi_{v^{(\ell-1)}}=0$, and (\ref{eq:bra}) and (\ref{3-F}) yield
$$
\EE_2=\partialx{v^{(\ell-2)}}\;,\ \ \ 
[\XX,\EE_2]=-\,\alpha_{v^{(\ell-2)}}\,\partialx{z}\ .
$$
Since $[\XX,\EE_2]=\mu\XX_1$ (see (\ref{bralie2})), relations (\ref{cro1}) and
(\ref{taux3}) imply that
$\alpha_{v^{(\ell-2)}}$, $\mu$, and the bracket $[\XX,\EE_2]$ are zero, and
prove the first part of the lemma. Let us turn to $[\XX,\EE_3]$~:
from (\ref{eq:bra}) and (\ref{3-F}), one gets, since $\EE_2$ and $\XX$
commute, and $\XX\chi=0$, 
 \begin{equation}
   \label{eq:EE3}
   \EE_3=\chi_{v^{(\ell-2)}}\partialx{u^{(k-3)}}+\partialx{v^{(\ell-3)}}
\;,\ \ \ 
[\XX,\EE_3]=-(\EE_3\alpha)\,\partialx{z}\ .
 \end{equation}

In order to prove that $\EE_3\alpha=0$, 
let us examine equation (\ref{det0}). For short, we use the symbol $\mathcal{O}$
to denote \emph{any function} that depends on
$u,\ldots,u^{(k-3)},x,y,z,v,\ldots,v^{(\ell-3)}$ only. For instance,
$\XX\gamma-\delta=\mathcal{O}$, and all terms in the determinant are of this
nature, except the following three~:
\begin{eqnarray*}
  \FF\alpha & = & \chi\,\alpha_{u^{(k-3)}} +v^{(\ell-2)}\alpha_{v^{(\ell-3)}} +\mathcal{O},\\
\XX\FF\alpha & = & \chi\,\XX\alpha_{u^{(k-3)}} +v^{(\ell-2)}\XX\alpha_{v^{(\ell-3)}}
+\mathcal{O},\\
\XX_1\alpha & = & -\alpha_z\,\left(\chi\,\alpha_{u^{(k-3)}}
  +v^{(\ell-2)}\alpha_{v^{(\ell-3)}}\right)  +\mathcal{O}
\end{eqnarray*}
(we used $\XX\chi=0$).
Setting $\zeta=\chi\,\alpha_{u^{(k-3)}}
+v^{(\ell-2)}\alpha_{v^{(\ell-3)}}$, one has
\begin{equation}
  \label{99}
  \XX\zeta =\frac{\XX\alpha_{u^{(k-3)}}}{\alpha_{u^{(k-3)}}}\zeta
+\mathbf{b}\,v^{(\ell-2)}
\ \ \ \mbox{with}\ \ \ 
\mathbf{b}=\XX\alpha_{v^{(\ell-3)}}
-\alpha_{v^{(\ell-3)}}\,\frac{\XX\alpha_{u^{(k-3)}}}{\alpha_{u^{(k-3)}}}\ ,
\end{equation}
and equation (\ref{det0}) reads
\begin{eqnarray}
\label{88}
  \zeta^2
+\mathcal{O}\,\zeta
-(\XX\gamma-\delta)\,\mathbf{b}\,v^{(\ell-2)}+\mathcal{O}
&=&0\ .
\end{eqnarray}
Differentiating with respect to $\XX$ and using (\ref{99}) yields
\begin{eqnarray*}
2\,  \frac{\XX\alpha_{u^{(k-3)}}}{\alpha_{u^{(k-3)}}}\,\zeta^2
+
\left(2\mathbf{b}\,v^{(\ell-2)}+
\mathcal{O}\right)
\zeta
+\mathcal{O}\, 
v^{(\ell-2)}
+\mathcal{O} &=&0\ .
\end{eqnarray*}
Then, eliminating $\zeta$ between these two polynomials yields the resultant
$$
\left|
  \begin{array}{cccc}
1 & \mathcal{O} & -(\XX\gamma-\delta)\mathbf{b}\, v^{(\ell-2)}+\mathcal{O} & 0
\\
0 & 1 & \mathcal{O} & -(\XX\gamma-\delta)\mathbf{b}\, v^{(\ell-2)}+\mathcal{O}
\\
2\frac{\XX\alpha_{u^{(k-3)}}}{\alpha_{u^{(k-3)}}} & 
2\mathbf{b}\, v^{(\ell-2)}+\mathcal{O} &
\mathcal{O}\, v^{(\ell-2)}+\mathcal{O} & 0
\\
0 & 2\frac{\XX\alpha_{u^{(k-3)}}}{\alpha_{u^{(k-3)}}} & 
2\mathbf{b}\, v^{(\ell-2)}+\mathcal{O} &
\mathcal{O}\, v^{(\ell-2)}+\mathcal{O}
  \end{array}
\right|=0\ .
$$
This is a polynomial of degree at most three with respect to $v^{(\ell-2)}$,
the coefficient of $(v^{(\ell-2)})^3$ being $-4\mathbf{b}^3(\XX\gamma-\delta)$.
Hence $\mathbf{b}=0$ and, from (\ref{88}), $\zeta$ does not depend on
$v^{(\ell-2)}$.
This implies $\EE_3\alpha=0$ because, from 
(\ref{eq:EE3}) and the definition of $\zeta$, one has
$
\zeta_{v^{(\ell-2)}}=\EE_3\alpha
$.
\end{proof}

\begin{proof}[Proof of Lemma \Rref{lem-3}]
The independent variables in $\sysedpl\gamma\delta33$ are $u,\dot{u},\ddot{u},x,v,\dot{v},\ddot{v}$.
Since the determinant (\ref{eq:det2}) is nonzero,
one defines local coordinates 
$(x,y,z,w,v,\dot{v},\ddot{v})$ by
\begin{equation}
  \label{yzw}
y=p(u,\dot{u},\ddot{u},x,v,\dot{v},\ddot{v}),\ \ \ 
z=p_x(u,\dot{u},\ddot{u},x,v,\dot{v},\ddot{v}),\ \ \ 
w=p_{xx}(u,\dot{u},\ddot{u},x,v,\dot{v},\ddot{v}).
\end{equation}
In these coordinates, $\XX$ and $\FF$, defined in (\ref{XY}), have the following expressions, with
$\psi$ and $\alpha$ some functions to be studied further~:
\begin{eqnarray}
  \label{4-X}
  \XX&=&\partialx{x}+z\partialx{y}+w\partialx{z}+\alpha\partialx{w}\;,
\\
  \label{4-F}
  \FF&=&\gamma\partialx{y}+\delta\partialx{z}
+\psi\partialx{w}
+\dot{v}\partialx{v}+\ddot{v}\partialx{\dot{v}}\;.
\end{eqnarray}
Then, using, for short, the following notation $\Gamma$~:
\begin{equation}
  \label{Gamma}
  \Gamma\ =\ \XX\gamma-\delta\ \neq\ 0\ ,
\end{equation}
one has
\begin{eqnarray}
\label{C4-X1}
\XX_1&\!\!\!\!=\!\!\!\!&\Gamma\partialx{y}+\left(\XX\delta-\psi\right)\partialx{z}
+\left(\XX\psi-\FF \alpha\right)\partialx{w},
\\
\label{C4-[X,X1]}
\!\!\!\!\!\!\!\!\!\!\!\!\!
[\XX,\XX_1]&\!\!\!\!=\!\!\!\!&\left(\XX\Gamma-\XX\delta+\psi\right)\partialx{y}
+\left(\XX^2\delta-2\XX\psi+\FF \alpha\right)\partialx{z}
+\left(\XX^2\psi-\XX\FF \alpha -\XX_1\alpha\right)\partialx{w}.
\end{eqnarray}

Also,
$$
\EE=\partialx{\ddot{v}}
\,,\ \ \ \ 
\EE_2=[\EE_1,\FF]=\psi_{\ddot{v}}\partialx{w}+\partialx{\dot{v}}
\,,\ \ \ \ 
[\XX,\EE_2]=
\psi_{\ddot{v}}\partialx{z}+\left(\XX\psi_{u^{(k-1)}}-\EE_2\alpha\right)\partialx{w}\,
$$
but, from Lemma~\ref{lem-2cro}, one has $[\XX,\EE_2]=0$, hence
$\psi_{\ddot{v}}=0$, $\EE_2=\partial/\partial{\dot{v}}$ and
$\alpha_{\dot{v}}=0$. 
Then
$$
\EE_3=[\partialx{\dot{v}},\FF]=\psi_{\dot{v}}\partialx{w}+\partialx{v}
\,,\ \ \ \ 
[\XX,\EE_3]=\psi_{\dot{v}}\partialx{z}+
\left(\XX\psi_{\dot{v}}-\EE_3\alpha\right)\partialx{w}\,,
$$
but, from Lemma~\ref{lem-2cro}, one has $[\XX,\EE_3]=0$, hence
$\psi_{\dot{v}}=0$, $\EE_3=\partial/\partial{v}$ and
$\alpha_{v}=0$. 
To sum up, 
\begin{equation}
  \label{eq:0}
\EE=\partialx{\ddot{v}}
\,,\ \ \ \ 
  \EE_2=\partialx{\dot{v}}\,,\ \ \EE_3=\partialx{v}\,,
\end{equation}
$\alpha$ depends at most on $(x,y,z,w)$ only and $\psi$ on $(x,y,z,w,v)$.
\\\textbf{Notation:}
until the end of this proof, $\mathcal{O}$ stands for \emph{any} function of
$x,y,z,w$ only. For instance, $\alpha=\mathcal{O}$, $\gamma=\mathcal{O}$,
$\delta=\mathcal{O}$, $\Gamma=\mathcal{O}$, $\XX\Gamma=\mathcal{O}$ ,
$\XX\delta=\mathcal{O}$ and $\XX^2\delta=\mathcal{O}$. 

From (\ref{bralie}), (\ref{C4-X1}) and (\ref{C4-[X,X1]}), one has
$\left|
\begin{array}{cc}
\Gamma & \XX\delta-\psi\\
\XX\Gamma-\XX\delta+\psi & \XX^2\delta-2\XX\psi+\FF \alpha 
\end{array}\right|=0\ .$
Hence
\begin{equation}
  \label{Xpsi}
  \XX\psi\ =\ \frac1{2\Gamma}\psi^2+\mathcal{O}\psi+\mathcal{O}\ .
\end{equation}
We now write the expression (\ref{C4-X1}) of $\XX_1$ as
\begin{eqnarray}
\label{X1sep}
&&\XX_1=\XX_1^0+\psi\XX_1^1+\psi^2\XX_1^2
\\
  \label{X1011}
&\mbox{with}&\displaystyle
  \XX_1^0=\Gamma\partialx{y}+\mathcal{O}\,\partialx{z}+\mathcal{O}\,\partialx{w}\;,\ \ \ 
\XX_1^1=-\,\partialx{z}+\mathcal{O}\,\partialx{w}\;,\ \ \ 
\XX_1^2=\frac1{2\Gamma}\partialx{w}\;.
\end{eqnarray}
Note that $\XX_1^0$, $\XX_1^1$ and $\XX_1^2$ are vector fields in the variables $x,y,z,w$ only.
Now define~:
\begin{eqnarray}
  \label{U}
  U&=& -\XX_1^1-\frac{\psi}{\Gamma}\,\partialx{w}\ =\
  \partialx{z}+\left(\mathcal{O}-\frac\psi\Gamma\right)\partialx{w}\ ,
\\
  \label{V}
  V&=& \XX_1^0-\psi^2\XX_1^2\ =\ 
\Gamma\partialx{y}+\mathcal{O}\,\partialx{z}+\left(\mathcal{O}-\frac{\psi^2}{2\Gamma}\right)\partialx{w}\
,
\end{eqnarray}
so that
\begin{equation}
  \label{X1UV}
  \XX_1=V-\psi U
\end{equation}
and, from (\ref{4-F}) and (\ref{X1sep}) one deduces the following expression of $\XX_2=[\XX_1,\FF]$~:
\begin{equation}
\label{C4-X2}
\XX_2\ =\ \big(\FF\psi\big)\,U+\big(\XX_1\psi\big)\,\partialx{w}
+\psi^3\frac{\Gamma_w}{2\Gamma^2}\,\partialx{w}
+\psi^2\left(\frac{\gamma_w}{2\Gamma}\partialx{y}+\mathcal{O}\partialx{z}+\mathcal{O}\partialx{w}\right)
+\psi\XX_2^1+\XX_2^0
\end{equation}
where $\XX_2^1$ and $\XX_2^0$ are two vector fields in the variables
$x,y,z,w$ only.

This formula and (\ref{eq:0}) imply
$[\EE_2,\XX_2]=\big(\FF\psi\big)_{\dot{v}}\,U=\psi_{v}\,U$; hence, from
the last relation in (\ref{bralie2}), either $\psi_{v}$ is identically zero or
$U$ is a linear combination of $\XX_1$ and $\XX_2$. 
We assume, until the end of the proof, that $U$ is a linear combination of
$\XX_1$ and $\XX_2$. 
This implies, using (\ref{X1UV}), that $\XX_2$ and $\XX_1$ are linear
combinations of $U$ and $V$; hence $U,V$ is
another basis for $\XX_1,\XX_2$. 
Also, from (\ref{bralie}) $[U,V]$ must be a linear combination of $U$ and $V$.
From (\ref{U}) and (\ref{V}),
\begin{displaymath}
  [U,V]=\frac{\XX_1\psi}{\Gamma}\partialx{w}
-\psi^2\,\mathcal{O}\,\partialx{w}
+\psi\,W^1+W^0
\end{displaymath}
where $W^1$ and $W^0$ are two vector fields in the variables
$x,y,z,w$ only, and, finally, with $Z^1$ and $Z^0$ two other vector fields in the variables
$x,y,z,w$ only, one has, from (\ref{C4-X2})
\begin{displaymath}
  \XX_2-(\FF\!\psi)\,U-\Gamma\,[U,V]\ =\ 
\psi^3\frac{\Gamma_w}{2\Gamma^2}\,\partialx{w}
+
\psi^2\left(\frac{\gamma_w}{2\Gamma}\partialx{y}
+\mathcal{O}\partialx{z}+\mathcal{O}\partialx{w}\right)
+\psi\,Z^1+Z^0\ .
\end{displaymath}
This vector field is also a linear
combination of $U$ and $V$. Computing the determinant in the basis
$\partial/\partial y$, $\partial/\partial z$, $\partial/\partial w$, one has,
using (\ref{U}) and (\ref{V}),
$$
\det\big(U,V,\,\XX_2-(\FF\!\psi)U-\Gamma[U,V]\,\big)=
\frac{\gamma_w}{\Gamma^3}\psi^4+\mathcal{O}\psi^3+\mathcal{O}\psi^2+\mathcal{O}\psi+\mathcal{O}
=0\
.
$$
It is assumed from the definition of $\sysedp$ that the
partial derivative of $\gamma$ with respect to its fourth argument is
nonzero; hence $\gamma_w\neq0$ and the above polynomial of degree 4 with
respect to $\psi$ is nontrivial; its coefficients depend on $x,y,z,w$ only,
hence $\psi$ cannot depend on $v$. 

We have proved that, in any case, both $\alpha$ and $\psi$ depend on
$x,y,z,w$ only, and this yields the desired identities in the lemma.
\end{proof}

\section{}
\label{app-lem0}
\begin{lmm}
\label{lem0}
Let $p$ be a smooth function of $u,\ldots,u^{(k-1)},x,v,\ldots,v^{(\ell-1)}$, $r$
a smooth function of
$u,\ldots,u^{(k-1)},x$, $v,\ldots,v^{(\ell)}$, with $r_{v^{(\ell)}}\neq0$, and $f$ 
a smooth function of four variables such that
\begin{equation}
  \label{ducon0}
  \sum_{i=0}^{k-2}u^{(i+1)}p_{u^{(i)}}\;+\;rp_{u^{(k-1)}}\;+\;
  \sum_{i=0}^{\ell-1}v^{(i+1)}p_{v^{(i)}}\ =\ f(x,p,p_{x},p_{xx})
\end{equation}
where, by convention, $rp_{u^{(k-1)}}$ is zero if $k=0$ and the first
(resp. last) sum is zero if $k\leq 1$ (resp. $\ell=0$).
Then either $p$ depends on $x$ only or
\begin{equation}
  \label{ducon}
  k\geq1\,,\ \ \ \ell\geq1\,,\ \ \ p_{u^{(k-1)}}\neq0\,,\ \ \ p_{v^{(\ell-1)}}\neq0\,.
\end{equation}
\end{lmm}
\begin{proof}
Let $m\leq k-1$ and $n\leq\ell-1$ be the smallest integers such that
$p$ depends on $u,\ldots,u^{(m)},x,v,\ldots,v^{(n)}$;
if $p$ depends on none of the variables $u,\ldots,u^{(k-1)}$ (or
$v,\ldots,v^{(\ell-1)}$), take $m<0$ (or $n<0$).
Then $p_{u^{(m)}}\neq0$ if $m\geq0$
and $p_{v^{(n)}}\neq0$ if $n\geq0$.

The lemma states that either $m<0$ and $n<0$ or
$k\geq 1$, $\ell\geq1$ and $(m,n)=(k-1,\ell-1)$.
This is indeed true~:\\
- if $m=k-1$ and $k\geq1$ then
$n=\ell-1$ and $\ell\geq 1$ because if not, differentiating both sides in (\ref{ducon0})
with respect to $v^{(\ell)}$ would yield $r_{v^{(\ell)}}p_{u^{(k-1)}}=0$, but
the lemma assumes that $r_{v^{(\ell)}}\neq0$,\\
- if $m<k-1$ or $m=0$, (\ref{ducon0}) becomes~:
$\sum_{i=0}^{m}u^{(i+1)}p_{u^{(i)}}+\sum_{i=0}^{n}v^{(i+1)}p_{v^{(i)}}=f(x,p,p_{x},p_{xx})$;
if $m\geq0$, differentiating with respect to $u^{(m+1)}$ yields $p_{u^{(m)}}=0$ and if 
$n\geq0$, differentiating with respect to $u^{(m+1)}$ yields $p_{v^{(n)}}=0$;
hence $m$ and $n$ must both be negative. 
\end{proof}

\bigskip

\bigskip
\noindent
This paper owes a lot not only to the original article
\cite{Hilb12}, but also to the careful re-reading  of that article by \emph{P.~Rouchon} in
\cite{Rouc92} (see also in \cite{Rouc94}).

\medskip

\noindent
  Also, the authors are very grateful to their colleague \emph{José Grimm} at INRIA Sophia
  Antipolis for an extremely careful reading of the manuscript that led to many
  improvements.


\begin{thebibliography}{10}

\bibitem{Aran-Moo-Pom95vars}
E.~Aranda-Bricaire, C.~H. Moog, and J.-B. Pomet.
\newblock An infinitesimal {B}runovsky form for nonlinear systems with
  applications to dynamic linearization.
\newblock {\em Banach Center Publications}, 32:19--33, 1995.

\bibitem{Avan05th}
D.~Avanessoff.
\newblock {\em Linéarisation dynamique des systèmes non linéaires et
  paramétrage de l'ensemble des solutions}.
\newblock PhD Thesis, Univ. de Nice - Sophia Antipolis, June 2005.

\bibitem{Brya-Che-Gar-G-G91}
R.~L. Bryant, S.~S. Chern, R.~B. Gardner, H.~L. Goldschmitt, and P.~A.
  Griffiths.
\newblock {\em Exterior Differential Systems}, volume~18 of {\em Mathematical
  Sciences Research Institute Publications}.
\newblock Springer-Verlag, 1991.

\bibitem{Cart15}
{\'E}.~Cartan.
\newblock Sur l'int{\'e}gration de certains syst{\`e}mes ind{\'e}termin{\'e}s
  d'{\'e}quations diff{\'e}rentielles.
\newblock {\em J. f{\"u}r reine und angew. Math.}, 145:86--91, 1915.

\bibitem{Char-Lev-Mar89}
B.~Charlet, J.~L{\'e}vine, and R.~Marino.
\newblock On dynamic feedback linearization.
\newblock {\em Syst. \& Control Lett.}, 13:143--151, 1989.

\bibitem{Char-Lev-Mar91}
B.~Charlet, J.~L{\'e}vine, and R.~Marino.
\newblock Sufficient conditions for dynamic state feedback linearization.
\newblock {\em SIAM J. on Control and Optim.}, 29:38--57, 1991.

\bibitem{Flie-Lev-Mar-R92cras}
M.~Fliess, J.~L{\'e}vine, P.~Martin, and P.~Rouchon.
\newblock Sur les syst{\`e}mes non lin{\'e}aires diff{\'e}rentiellement plats.
\newblock {\em C.~R. Acad. Sci. Paris}, S{\'e}rie I, 315:619--624, 1992.

\bibitem{Flie-Lev-Mar-R95ijc}
M.~Fliess, J.~L{\'e}vine, P.~Martin, and P.~Rouchon.
\newblock Flatness and defect of nonlinear systems: Introductory theory and
  examples.
\newblock {\em Int. J. Control}, 61(6):1327--1361, 1995.

\bibitem{Flie-Lev-Mar-R99geo}
M.~Fliess, J.~L{\'e}vine, P.~Martin, and P.~Rouchon.
\newblock A {L}ie-{B}\"acklund approach to equivalence and flatness of
  nonlinear systems.
\newblock {\em IEEE Trans. Automat. Control}, 44(5):922--937, 1999.

\bibitem{Flie-Lev-Mar-R99open}
M.~Fliess, J.~L{\'e}vine, P.~Martin, and P.~Rouchon.
\newblock Some open questions related to flat nonlinear systems.
\newblock In {\em Open problems in mathematical systems and control theory},
  pp. 99--103. Springer, London, 1999.

\bibitem{Golu-Gui73}
M.~Golubitsky and V.~Guillemin.
\newblock {\em Stable mappings and their singularities}.
\newblock Springer-Verlag, New York, 1973.
\newblock GTM, Vol. 14.

\bibitem{Hilb12}
D.~Hilbert.
\newblock {\"U}ber den {B}egriff der {K}lasse von {D}ifferentialgleichungen.
\newblock {\em Math. Annalen}, 73:95--108, 1912.

\bibitem{Hube03}
E.~Hubert.
\newblock Notes on triangular sets and triangulation-decomposition algorithms.
  {I}: {P}olynomial systems. {II}: {D}ifferential systems.
\newblock In F.~Winkler et al. eds., {\em Symbolic and Numerical
  Scientific Computing}, LNCS Vol. 2630,
  pp. 1--87. Springer Verlag, 2003.

\bibitem{Isid-Moo-deL86}
A.~Isidori, C.~H. Moog, and A.~de~Luca.
\newblock {A} sufficient condition for full linearization via dynamic state
  feedback.
\newblock In {\em Proc. 25th IEEE Conf. on Decision \& Control, Athens}, pp.
  203--207, 1986.

\bibitem{Mart92th}
P.~Martin.
\newblock {\em Contribution {\`a} l'{\'e}tude des syst{\`e}mes
  differentiellement plats}.
\newblock PhD thesis, Ecole des Mines, Paris, 1992.

\bibitem{Mart-Rou-Mur02tri}
P.~Martin, R.~M. Murray, and P.~Rouchon.
\newblock Flat systems.
\newblock In {\em Mathematical control theory, Part 1, 2 (Trieste, 2001)}, ICTP
  Lect. Notes, VIII, pp. 705--768 (electronic). Abdus Salam Int. Cent.
  Theoret. Phys., Trieste, 2002.

\bibitem{Mart-Rou94}
P.~Martin and P.~Rouchon.
\newblock Feedback linearization and driftless systems.
\newblock {\em Math. of Control, Signals \& Systems}, 7:235--254, 1994.

\bibitem{Pome95vars}
J.-B. Pomet.
\newblock A differential geometric setting for dynamic equivalence and dynamic
  linearization.
\newblock {\em Banach Center Publications}, 32:319--339, 1995.

\bibitem{Pome97cocv}
J.-B. Pomet.
\newblock On dynamic feedback linearization of four-dimensional affine control
  systems with two inputs.
\newblock {\em ESAIM Control Optim. Calc. Var.} (this journal), 2:151--230, June 1997.
\newblock URL: http://www.edpsciences.org/cocv/.

\bibitem{Ritt50}
J.~F. Ritt.
\newblock {\em Differential {A}lgebra}.
\newblock AMS Colloq. Publ., Vol. XXXIII.
  New York, 1950.

\bibitem{Rouc92}
P.~Rouchon.
\newblock Flatness and oscillatory control: some theoretical results and case
  studies.
\newblock Tech. report PR412, CAS, Ecole des Mines, Paris, Sept. 1992.

\bibitem{Rouc94}
P.~Rouchon.
\newblock Necessary condition and genericity of dynamic feedback linearization.
\newblock {\em J. of Math. Systems, Estimation, and Control}, 4:1--14, 1994.

\bibitem{Slui93}
W.~M. Sluis.
\newblock A necessary condition for dynamic feedback linearization.
\newblock {\em Syst. \& Control Lett.}, 21:277--283, 1993.

\bibitem{VanN-Rat-Mur98}
M.~van Nieuwstadt, M.~Rathinam, and R.~Murray.
\newblock Differential flatness and absolute equivalence of nonlinear control
  systems.
\newblock {\em SIAM J. on Control and Optim.}, 36(4):1225--1239, 1998.
\newblock http://epubs.siam.org:80/sam-bin/dbq/article/27402.

\bibitem{Zerv32}
P.~Zervos.
\newblock Le problème de {M}onge.
\newblock {\em Mémorial des Sciences Mathématiques}, LIII, 1932.

\end{thebibliography}
\end{document}